\documentclass[a4paper]{article}
\usepackage{graphicx}
\usepackage[utf8]{inputenc}
\usepackage[ngerman, english]{babel}
\usepackage{amsmath,amsthm,amssymb,amsfonts}
\usepackage{mathrsfs}
\usepackage{ifthen}
\usepackage{enumerate}
\usepackage{comment}
\usepackage{subcaption}
\usepackage{hyphsubst}
\usepackage{xcolor}
\usepackage{calc}
\usepackage{graphicx}
\usepackage[toc,page]{appendix}

\makeatletter
\def\@hspace#1{\begingroup\setlength\dimen@{#1}\hskip\dimen@\endgroup}
\makeatother

\newtheorem{theorem}{Theorem}[section]
\newtheorem{lemma}[theorem]{Lemma}
\newtheorem{definition}[theorem]{Definition}

\newtheorem{remark}[theorem]{Remark}
\newtheorem{corollary}[theorem]{Corollary}
\newtheorem{assumption}[theorem]{Assumption}
\newtheorem{proposition}[theorem]{Proposition}

\newcommand{\signv}{\fraks_V}
\newcommand{\signw}{\fraks_{W_0}}
\newcommand{\poop}{\pi}
\newcommand{\ipop}{\Pi}
\newcommand{\hn}{h_{n}}
\newcommand{\CR}{\rho}
\newcommand{\CC}{C_\ipop}

\newcommand{\nv}{\nu}
\newcommand{\tv}{\tau}
\newcommand{\ol}[1]{\overline{#1}}

\newcommand{\spl}{\langle}
\newcommand{\spr}{\rangle}
\newcommand{\bpm}{\begin{pmatrix}}
\newcommand{\epm}{\end{pmatrix}}

\DeclareMathOperator{\curl}{curl}
\DeclareMathOperator{\Curl}{Curl}

\renewcommand{\div}{\operatorname{div}}

\DeclareMathOperator{\ind}{ind}

\DeclareMathOperator{\supp}{supp}

\newcommand{\setC}{\mathbb{C}}

\newcommand{\setN}{\mathbb{N}}

\newcommand{\setQ}{\mathbb{Q}}
\newcommand{\setR}{\mathbb{R}}

\newcommand{\calP}{\mathcal{P}}

\newcommand{\calR}{\mathcal{R}}

\newcommand{\calT}{\mathcal{T}}

\newcommand{\fraks}{\mathfrak{s}}

\definecolor{brickred}{rgb}{0.8, 0.25, 0.33}
\definecolor{bostonuniversityred}{rgb}{0.8, 0.0, 0.0}
\definecolor{cornellred}{rgb}{0.7, 0.11, 0.11}
\definecolor{corn}{rgb}{0.98, 0.93, 0.36}
\definecolor{schoolbusyellow}{rgb}{1.0, 0.85, 0.0}
\definecolor{TUblue}{rgb}{0,102,153}
\colorlet{TUbluelight}{TUblue!30!white}

\title{On the approximation of dispersive electromagnetic eigenvalue problems in 2D}
\author{Martin Halla\footnote{Max-Planck-Institut f\"ur Sonnensystemforschung, Justus-von-Liebig-Weg 3, 37077 G\"ottingen, Deutschland (halla@mps.mpg.de)}
\footnote{Institut f\"ur Numerische und Angewandte Mathematik, Georg-August Universit\"at G\"ottingen, Lotzestra\ss e 16-18, 37083 Göttingen, Deutschland}}
\date{\today}

\begin{document}
\maketitle

\begin{abstract}
\noindent
We consider time-harmonic electromagnetic wave equations in composites of a dispersive material surrounded by a classical material.
In certain frequency ranges this leads to sign-changing permittivity and/or permeability.
Previously meshing rules were reported, which guarantee the convergence of finite element approximations to the related scalar source problems.
Here we generalize these results to the electromagnetic two dimensional vectorial equations and the related holomorphic eigenvalue problems.
Different than for the analysis on the continuous level, we require an assumption on both contrasts of the permittivity and the permeability.
We confirm our theoretical results with computational studies.\\

\noindent
Keywords: Maxwell’s equations; eigenvalues; metamaterials; plasmonics; sign-changing coefficients; compatible discretization.
\end{abstract}

\section{Introduction}\label{sec:introduction}

For dispersive material laws \cite{Pendry:00,ZiolkowskiHeyman:01,GralakTip:10,CuiSmithLiu:10} the magnetic permeability and/or the electric permittivity are negative in certain frequency ranges.
In composites of classical and dispersive materials this leads to indefinite coefficients for the time-harmonic wave equations
\begin{align*}
-\div(\mu(\omega)^{-1} \nabla u)-\omega^2\epsilon(\omega)u=0, \qquad
\curl(\mu(\omega)^{-1} \curl u)-\omega^2\epsilon(\omega)u=0
\end{align*}
and the well-posedness/Fredholmness of the associated boundary value problems/operators becomes a delicate question which depends on the contrasts $\kappa_{\mu^{-1}}$, $\kappa_\epsilon$ at the materials interface.
In \cite{BonnetBDChesnelCiarlet:12} the Fredholmness for the scalar wave equation was thoroughly answered.
For smooth interfaces the operator is Fredholm if $\kappa_{\mu^{-1}}\neq1$.
If the interface admits corners, then the operator is Fredholm if the contrast $\kappa_{\mu^{-1}}$ is outside a critical interval $I_c\ni-1$, whose size depends on the angles of the corners.
We refer to \cite{BonnetBDCarvalhoChesnelCiarlet:15} for more details in the presence of corners and a technique to restore a Fredholm framework for contrasts inside of the critical interval.
In \cite{BonnetBDChesnelCiarlet:14b,BonnetBDChesnelCiarlet:14a} the analysis of \cite{BonnetBDChesnelCiarlet:12} was extended to the two- and three-dimensional electromagnetic equations.
Instead of the geometric techniques used in \cite{BonnetBDChesnelCiarlet:12}, the technique of \cite{BonnetBDChesnelCiarlet:14b,BonnetBDChesnelCiarlet:14a} was to relate the Fredholmness to the well-posedness of the associated scalar equations (already discussed in \cite{BonnetBDChesnelCiarlet:12}).
The finite element approximation of scalar source problems in $\setR^2$ was discussed in \cite{BonnetBDCarvalhoCiarlet:18} wherein meshing rules were reported which guarantee the convergence of approximations.

The applied technique of the former references to prove that an operator $A$ is Fredholm, is to construct a bijective operator $T$ such that $A$ is weakly $T$-coercive.
To explain this terminology recall that an operator $A\in L(X)$ is called coercive, if $\inf_{u\in X,\|u\|_X=1}|\spl Au,u\spr_X|>0$ and weakly coercive if there exists a compact operator $K\in L(X)$ such that $A+K$ is coercive.
For a bijective operator $T\in L(X)$ the operator $A\in L(X)$ is called (weakly) $T$-coercive, if $T^*A$ is (weakly) coercive.
Here $T^*$ denotes the adjoint operator of $T$.
Further we call Galerkin approximations with finite dimensional spaces $X_n\subset X$ to be $T$-compatible, if there exist operators $T_n\in L(X_n)$ with $\lim_{n\to\infty} \|T-T_n\|_{L(X_n,X)}=0$.
This criterion does not only guarantee the convergence of approximations of source problems $Au=f$, but also of holomorphic eigenvalue problems $A(\omega)u=0$, see the framework \cite{Halla:19Diss,Halla:19Tcomp} which is based on \cite{Karma:96a,Karma:96b}.
We also refer to the introduction of \cite{Halla:19Tcomp} for a historical recap of \emph{T-analysis}.

In this article we consider the the two dimensional vectorial electromagnetic wave equation and the finite element approximation of the associated source and holomorphic eigenvalue problems.
Different than for the analysis on the continuous level in 2D, we require an assumption on both contrasts of the permittivity and the permeability.
We prove that the meshing rules of \cite{BonnetBDCarvalhoCiarlet:18} guarantee the convergence of h-finite element approximations with N\'ed\'elec elements.
To this end we generalize the approach of \cite{BonnetBDCarvalhoCiarlet:18} and construct an apt $T$-operator by means of local pattern based reflection operators.
In particular we employ a Helmholtz decomposition of $H_0(\curl)=V\oplus^\bot \nabla H^1_0$ and introduce individual local reflection operators for each Helmholtz component.

The remainder of this article is structured as follows.
In Section~\ref{sec:setting} we specify our setting and introduce the operator $A$ to be investigated.
In Section~\ref{sec:wTc} we construct an apt $T$-operator and prove weak $T$-coercivity of $A$.
In Section~\ref{sec:Tcomp} we construct apt approximations $T_n$ of $T$ and prove the convergence of approximations.
In Section~\ref{sec:compexamp} we conduct computational studies to verify our theoretical results.
In Section~\ref{sec:conclusion} we conclude and discuss the further outlook.

\section{Specification of the problem}\label{sec:setting}

Let $\Omega\subset\setR^2$ be a bounded Lipschitz domain (open and connected) with unit outward normal vector $\nv$ and unit tangential vector $\tv$.
Let $\Omega_-\subset\Omega$ be a Lipschitz polyhedron with $\partial\Omega_-\cap\partial\Omega=\emptyset$ and let $\Omega_+:=\Omega\setminus\ol{\Omega_-}$, $\Sigma:=\partial\Omega_-$.
Let $\mu$ and $\epsilon$ be the magnetic permeability and the electric permittivity.
We consider the case of a composite material whereby $\Omega_+$ is a classical material and $\Omega_-$ is a dispersive material, i.e.\ $\epsilon|_{\Omega_-}$ and $\mu|_{\Omega_-}$ depend on the frequency $\omega$.
We formulate all our theory in terms of the spectral parameter $\lambda:=\omega^2$.
To simplify the presentation we assume that $\mu, \epsilon$ are scalar and constant in $\Omega_\pm$ and that in $\Omega_-$ the material is modeled by a Drude law.
To be precise let $\mu_+, \mu_-, \epsilon_+, \epsilon_-$ be positive constants and $\omega_\mu, \omega_\epsilon$ be real non-negative constants and
\begin{align*}
\mu(\lambda)|_{\Omega_+}&=\mu_+, \qquad
\mu(\lambda)|_{\Omega_-}=\mu_-\bigg(1-\frac{\omega_\mu^2}{\lambda}\bigg), \\
\epsilon(\lambda)|_{\Omega_+}&=\epsilon_+, \qquad
\epsilon(\lambda)|_{\Omega_-}=\epsilon_-\bigg(1-\frac{\omega_\epsilon^2}{\lambda}\bigg), \quad
\end{align*}
whereby $\omega_\mu$ and $\omega_\epsilon$ are the resonant frequencies.
We introduce the contrasts
\begin{align*}
\kappa_{\mu^{-1}}(\lambda):=\frac{\mu_-^{-1}}{\mu_+^{-1}}
\bigg(1-\frac{\omega_\mu^2}{\lambda}\bigg)^{-1}
\qquad\text{and}\qquad
\kappa_{\epsilon}(\lambda):=\frac{\epsilon_-}{\epsilon_+}
\bigg(1-\frac{\omega_\epsilon^2}{\lambda}\bigg).
\end{align*}
Note that our forthcoming theory in Sections~\ref{sec:wTc} and \ref{sec:Tcomp} can easily be adapted to non-homogeneous coefficients $\mu_\pm, \epsilon_\pm$, although in this case the definition of the contrasts is more technical (see, e.g.\ \cite[Theorem 4.3]{BonnetBDChesnelCiarlet:12}).
Let $\partial_{x_n}u$ be the partial derivative of a function $u$ with respect to the variable $x_n$.
For a two-vector function $u=(u_1,u_2)^\top$ let $\curl u:=\partial_{x_1} u_2-\partial_{x_2} u_1$.
For a scalar function $u$ let $\Curl u:=(\partial_{x_2}u,-\partial_{x_1}u)^\top$.
For $s\geq0$ and a Lipschitz domain $D\subset\Omega$ let $H^s(\curl;D):=\{u\in (H^s(D)^2\colon \curl u\in H^s(D)\}$ with scalar product $\spl u,u'\spr_{H^s(\curl;D)}:=\spl \curl u,\curl u' \spr_{H^s(D)}+\spl u, u' \spr_{H^s(D)}$, $u,u'\in H^s(\curl;D)$ and $H(\curl;D):=H^0(\curl;D)$.
We consider the electromagnetic time-harmonic wave equation
\begin{subequations}\label{eq:PDE}
\begin{align}
\Curl\mu(\lambda)^{-1}\curl u -\lambda\epsilon(\lambda) u &= f \quad\text{in }\Omega,\\
\tv\cdot u &= 0 \quad\text{on }\partial\Omega.
\end{align}
\end{subequations}
Here we consider either a source problem for which a solution $u\in H(\curl;\Omega)$ is sought for given $\lambda$ and right hand side $f$.
Or we consider the eigenvalue problem for which $f=0$ and a pair $(\lambda,u)$ is sought (with non-trivial $u$).
To formulate these in a functional setting we introduce the following.
For two Banach spaces $X,Y$ denote $L(X,Y)$ the space of bounded linear operators from $X$ to $Y$ with norm $\|A\|_{L(X,Y)}:=\sup_{u\in X,\|u\|_X=1}\|Au\|_Y$ and set $L(X):=L(X,X)$.
Further let
\begin{align*}
X:=H_0(\curl;\Omega):=\{u\in H(\curl;\Omega)\colon \,\tv\cdot u=0\text{ on }\partial\Omega\}
\end{align*}
with scalar product $\spl u,u'\spr_X:=\spl u,u'\spr_{H(\curl;\Omega)}$.
We consider the operator $A(\lambda)\in L(X)$ defined by
\begin{align}\label{eq:DefA}
\spl A(\lambda)u,u' \spr_X:=
\spl \mu(\lambda)^{-1} \curl u,\curl u' \spr_{L^2(\Omega)}
-\lambda \spl \epsilon(\lambda) u, u' \spr_{L^2(\Omega)}
\end{align}
for all $u,u'\in X$ and $\lambda\in\setC\setminus\{0,\omega_\mu^2\}$.
Then the operator formulation of the source problem is to find for given $\lambda\in\setC\setminus\{0,\omega_\mu^2\}$ and $\tilde f\in X$ a solution $u\in X$ such that $A(\lambda)u=\tilde f$.
The operator formulation of the eigenvalue problem is to find $(\lambda,u) \in (\setC\setminus\{0,\omega_\mu^2\})\times (X\setminus\{0\})$ such that $A(\lambda)u=0$.

In this article we are concerned with two goals.
The first is to prove that $A(\lambda)$ is Fredholm with index zero for $\lambda\in\Lambda$, whereby $\Lambda\subset\setC$ is an open connected set which will be specified later on.
The second one is two prove that $A(\lambda)$ can be sufficiently approximated by convenient $H(\curl)$-finite element methods with the meshing rules introduced in \cite{BonnetBDCarvalhoCiarlet:18}.
The means to establish the former goals are to construct for each $\lambda\in\Lambda$ a bijective operator $T(\lambda)\in L(X)$ such that $A(\lambda)$ is weakly $T(\lambda)$-coercive.
In this situation we say that $A(\cdot)$ is weakly $T(\cdot)$-coercive, from which it readily follows that $A(\lambda)$ is Fredholm with index zero for all $\lambda\in\Lambda$.
For the analysis of approximations of source and eigenvalue problems by Galerkin methods with discrete spaces $X_n\subset X$ we employ the framework of $T(\cdot)$-compatible approximations \cite{Halla:19Tcomp,Halla:19Diss}.
To this end we need to construct operators $T_n(\lambda)\in L(X_n)$ such that $\lim_{n\to\infty} \|T(\lambda)-T_n(\lambda)\|_{L(X_n,X)}=0$ for each $\lambda\in\Lambda$.

Let us examine the properties of $A(\lambda)$ for $\lambda$ in different regions of $\setC$.
At first we note the singular behaviour at $\lambda=0,\omega^2_\mu$ and the degenerate behaviour at $\lambda=\omega^2_\epsilon$.
On the other hand $A(\lambda)$ is clearly coercive for $\lambda<0$.
Also for $\lambda\in\setC\setminus\setR$ one can check that the coefficients
$1/\mu(\lambda)|_{\Omega_+}$, $-\lambda\epsilon(\lambda)|_{\Omega_+}$, $1/\mu(\lambda)|_{\Omega_-}$, $-\lambda\epsilon(\lambda)|_{\Omega_-}$ are contained in a closed salient sector.
Hence for $\lambda\in\setC\setminus\setR$ the operator $A(\lambda)$ is coercive too.
Thus with $\setR_+:=\{x\in\setR\colon x>0\}$ we set
\begin{align}\label{eq:DefT-identity}
T(\lambda):=I_X, \qquad \lambda\in\setC\setminus\setR_+.
\end{align}
For the analysis of $A(\lambda)$ for $\lambda>0$ we introduce the following Helmholtz decomposition.
Let
\begin{align*}
W_0&:=H^1_0(\Omega),\qquad
W:=\{\nabla w_0\colon w_0\in W_0\}\subset X,\\
V&:=W^{\bot_X}=\{u\in X\colon \div u=0 \text{ in }\Omega\},
\end{align*}
and $P_W, P_V$ be the respective orthogonal projections and $P_{W_0}\in L(X,W_0)$ be such that $P_W u=\nabla P_{W_0}u$ for all $u\in X$.
We recall the embeddings
\begin{align}\label{eq:Vembeddings}
V
\lhook\joinrel\xrightarrow{\text{bounded}} (H^{1/2}(\Omega))^2
\lhook\joinrel\xrightarrow{\text{compact}} (L^2(\Omega))^2
\end{align}
whereby the first embedding can be seen as in \cite{Costabel:90}.
Thence for $\lambda>\max(\omega^2_\mu,\omega^2_\epsilon)$ the operator $A(\lambda)$ is weakly $T(\lambda)$-coercive with
\begin{align}\label{eq:DefT-classical}
T(\lambda):=P_V-P_W, \qquad \lambda>\max(\omega^2_\mu,\omega^2_\epsilon).
\end{align}
On the contrary if $\lambda\in(0,\max(\omega^2_\mu,\omega^2_\epsilon))$ the study of $A(\omega)$ requires much more sophisticated techniques and will be dealt with in Section~\ref{sec:wTc}.

\section{Weak T-coercivity}\label{sec:wTc}

In this section we consider the case $\lambda\in(0,\max(\omega^2_\mu,\omega^2_\epsilon))$.
For contrasts $\kappa_{\mu^{-1}}(\lambda)$, $\kappa_{\epsilon}(\lambda)$ outside a critical interval (which will be specified) we construct an apt operator $T(\lambda)$, such that $A(\lambda)$ is weakly $T(\lambda)$-coercive.
The definition of $T(\lambda)$ involves cut-off functions $\chi_n$ and local reflection operators $R_n^{V,\pm}, R_n^{W_0,\pm}$.
In Subsection~\ref{subsec:wTc-abstract} we formulate assumptions on those, introduce $T(\lambda)$ under this pretext and prove the weak $T(\lambda)$-coercivity of $A(\lambda)$.
In Subsection~\ref{subsec:wTc-explicit} we explicitly construct cut-off functions and local reflection operators which satisfy the former requirements.
In Subsection~\ref{subsec:wTc-results} we conclude and formulate our results.

\subsection{Abstract framework}\label{subsec:wTc-abstract}

We follow \cite[p.\ 813-814]{BonnetBDCarvalhoCiarlet:18} and introduce some geometry based tools which will be necessary to construct $T(\lambda)$.
Recall that $\Sigma$ is a polygonal line without endpoints.
Let $N$ be the number of its corners,
$(c_n)_{n=1,\dots,N}$ be its corners,
$(\alpha_n)_{n=1,\dots,N}$ be the corner angles measured in $\Omega_+$,
and $(e_n)_{n=1,\dots,N}$ be its edges.

\begin{assumption}[cut-off functions]\label{ass:abstractCutoff}
Let $\chi_n\in C^\infty(\ol{\Omega},[0,1])$, $n=1,\dots,2N$ be smooth cut-off functions.
We associate the first $N$ functions to the corners of $\Sigma$ and the second $N$ functions to the edges of $\Sigma$.
In particular for each $n=1,\dots,N$ let $\chi_n$ have a support localized in a neighborhood of the interface, such that $\chi_n=1$ in a neighborhood of the corner $c_n$.
Further let $\chi_n, n=1,\dots,N$ have mutually disjoint supports.
For each $n=1,\dots,N$ let $\chi_{N+n}$ have a support localized in a neighborhood of the edge $e_n$, such that $\chi_n=0$ on $\Sigma\setminus e_n$.
Additionally let $\sum_{n=1}^{2N} \chi_n=1$ in a neighborhood of $\Sigma$.
\end{assumption}

For $\chi_n$ satisfying Assumption~\ref{ass:abstractCutoff} we further define $S_n:=\supp \chi_n$ as the support of the cut-off function $\chi_n$.
Next we introduce our assumptions on the local reflection operators.
We introduce separate reflection operators for each Helmholtz component.
To this end let
\begin{align*}
X_\pm&:=\{u|_{\Omega_\pm} \colon u\in X\}, \quad \spl u,u'\spr_{X_\pm}:=\spl u,u'\spr_{H(\curl;\Omega_\pm)},\\
W_{0,\pm}&:=\{w_0|_{\Omega_\pm} \colon w_0\in W_0\}, \quad \spl w_0,w_0'\spr_{W_{0,\pm}}:=\spl w_0,w_0'\spr_{H^1(\Omega_\pm)}.
\end{align*}
If necessary we shorten the notation and write for a function $u$ defined on $\Omega$ and an operator $R$ acting from $X_{\pm}$ or $W_{0,\pm}$ the short form $Ru$ instead of $Ru|_{\Omega_\pm}$.
We note that in the following assumption the reflection operators $R_n^{V,\pm}$ act from $X_\pm$ (and not only from restrictions of $V$, i.e.\ $V_\pm$).

\begin{assumption}[local reflection operators]\label{ass:abstractReflection}
For each $n=1,\dots,2N$ let $\Omega_{n}\subset\Omega$ be a Lipschitz domain with $\ol{S_n}\subset\Omega_n$ and $R_n^{V,\pm}\in L(X_{\pm},X_{\mp})$, $R_n^{W_0,\pm} \in L(W_{0,\pm},W_{0,\mp})$ satisfy the following matching conditions on the traces
\begin{align*}
\tv\cdot R_n^{V,\pm}u|_{\Sigma\cap \Omega_n} = \tv\cdot u|_{\Sigma\cap \Omega_n}, \qquad
(R_n^{W_{0,\pm}}w_0)|_{\Sigma\cap \Omega_n} = w_0|_{\Sigma\cap \Omega_n},
\end{align*}
for all $u\in X_\pm$, $w_0\in W_{0,\pm}$.
In addition let there exist $t\in(0,1/2)$ such that
\begin{align*}
R_n^{V,\pm} P_V \in L(X_\pm,H^t(\Omega_\mp\cap\Omega_n))
\end{align*}
for each $n=1,\dots,2N$.
\end{assumption}
For $\chi_n$ and $R_n^{V,\pm}, R_n^{W_0,\pm}$ satisfying Assumptions~\ref{ass:abstractCutoff} and \ref{ass:abstractReflection} respectively let
\begin{align*}
\|R_n^{V,\pm}\|&:=
\sup_{u\in X_\pm, \| \chi_n^{1/2}\curl u \|_{L^2(\Omega_\pm\cap S_n)}=1}
\| \chi_n^{1/2}\curl R_n^{V,\pm}u \|_{L^2(\Omega_\mp\cap S_n)},\\
\|R^{V,\pm}\|&:=\max_{n=1,\dots,2N} \|R_n^{V,\pm}\|,\\
I_{R_V}&:=[-\|R^{V,-}\|,-1/\|R^{V,+}\|],
\end{align*}
and
\begin{align*}
\|R_n^{W_0,\pm}\|&:=
\sup_{w_0\in W_{0,\pm}, \| \chi_n^{1/2}\nabla w_0 \|_{L^2(\Omega_\pm\cap S_n)}=1}
\| \chi_n^{1/2}\nabla R_n^{W_0,\pm}w_0 \|_{L^2(\Omega_\mp\cap S_n)},\\
\|R^{W_0,\pm}\|&:=\max_{n=1,\dots,2N} \|R_n^{W_0,\pm}\|,\\
I_{R_{W_0}}&:=[-\|R^{W_0,-}\|,-1/\|R^{W_0,+}\|].
\end{align*}
Subsequently we define
\begin{align*}
T^{V,\pm}u&:=\left\{\begin{array}{ll}
\phantom{-}u|_{\Omega_+}-(1\mp1)\sum_{n=1}^{2N}\chi_n R_n^{V,-}u|_{\Omega_-}, &\text{in }\Omega_+,\\
-u|_{\Omega_-}+(1\pm1)\sum_{n=1}^{2N}\chi_n R_n^{V,+}u|_{\Omega_+}, &\text{in }\Omega_-,
\end{array}\right.,\\
T^{W_0,\pm}w_0&:=\left\{\begin{array}{ll}
\phantom{-}w_0|_{\Omega_+}-(1\mp1)\sum_{n=1}^{2N}\chi_n R_n^{W_0,-}w_0|_{\Omega_-}, &\text{in }\Omega_+,\\
-w_0|_{\Omega_-}+(1\pm1)\sum_{n=1}^{2N}\chi_n R_n^{W_0,+}w_0|_{\Omega_+}, &\text{in }\Omega_-,
\end{array}\right..
\end{align*}
In addition we set $T^{V,0}:=I_X$ and $T^{W_0,0}:=I_{W_0}$.
It is straightforward to see that $T^{V,\signv}\in L(X)$, $T^{W_0,\signw}\in L(W_0)$, for each $\signv,\signw\in\{+,-,0\}$.
It is then natural to consider the composed operator
\begin{subequations}\label{eq:DefTprime}
\begin{align}
T'(\lambda):=T^{V,\signv} P_V - \nabla T^{W_0,\signw} P_{W_0}
\end{align}
with
\begin{align}
\begin{aligned}
\begin{array}{ll}
\signv=0, &\text{if } \kappa_{\mu^{-1}}(\lambda)>0,\\
\signv=-, &\text{if } \kappa_{\mu^{-1}}(\lambda)<-1,\\
\signv=+, &\text{if } \kappa_{\mu^{-1}}(\lambda)\in(-1,0),
\end{array}
\text{and}\quad
\begin{array}{ll}
\signw=0, &\text{if } \kappa_\epsilon(\lambda)>0,\\
\signw=-, &\text{if } \kappa_\epsilon(\lambda)<-1,\\
\signw=+, &\text{if } \kappa_\epsilon(\lambda)\in(-1,0),
\end{array}
\end{aligned}
\end{align}
\end{subequations}
and observe $T'(\lambda)\in L(X)$.
Now we are able to proof a similar result as in Lemma~2 of \cite{BonnetBDCarvalhoCiarlet:18}.

\begin{lemma}\label{lem:TprimA-wc}
Let $\chi_n$ and $R_n^{V,\pm}, R_n^{W_{0,\pm}}$, $n=1,\dots,2N$ satisfy Assumptions~\ref{ass:abstractCutoff} and \ref{ass:abstractReflection} respectively.
Let $A(\lambda)$ be defined as in \eqref{eq:DefA} and $T'(\lambda)$ be defined as in \eqref{eq:DefTprime}.
Let $\lambda\in(0,\max(\omega^2_\mu,\omega^2_\epsilon))$ be such that $\kappa_{\mu^{-1}}(\lambda) \notin I_{R_V}$ and $\kappa_{\epsilon}(\lambda) \notin I_{R_{W_0}}$.
Then $(T'(\lambda))^*A(\lambda)$ is weakly coercive.
\end{lemma}
\begin{proof}
We present only the case $\kappa_{\mu^{-1}}(\lambda),\kappa_\epsilon(\lambda)\in(-1,0)$.
The other cases can be treated analogously.
We abbreviate $v:=P_V u$, $w_0:=P_{W_0} u$, $v':=P_V u'$ and $w_0:=P_{W_0} u'$ and define $A_1, A_2, A_3 \in L(X)$ by
\begin{align*}
\spl A_1u,u \spr_X&:=
\mu^{-1}_+ \spl \curl v,\curl v' \spr_{L^2(\Omega_+)}
+\mu^{-1}_-|1-\omega^2_\mu/\lambda|^{-1} \spl \curl v,\curl v' \spr_{L^2(\Omega_-)}\\
&-2\mu^{-1}_-|1-\omega^2_\mu/\lambda|^{-1} \sum_{n=1}^{2N} \spl \curl v,\chi_n \curl R_n^{V,+}v' \spr_{L^2(\Omega_-\cap S_n)}\\
&+\spl v, v' \spr_{L^2(\Omega)},\\
\spl A_2u,u \spr_X&:=\epsilon_+ \spl \nabla w_0,\nabla w_0' \spr_{L^2(\Omega_+)}
+\epsilon_-|1-\omega^2_\epsilon/\lambda| \spl \nabla w_0,\nabla w_0' \spr_{L^2(\Omega_-)}\\
&-2\epsilon_-|1-\omega^2_\epsilon/\lambda| \sum_{n=1}^{2N} \spl \nabla w_0,\chi_n \nabla R_n^{W_0,+}w_0' \spr_{L^2(\Omega_+\cap S_n)},\\
\spl A_3u,u \spr_X&:=
-2\mu^{-1}_-|1-\omega^2_\mu/\lambda|^{-1} \sum_{n=1}^{2N} \spl \curl v,(\nabla^\top\chi_n) \cdot R_n^{V,+}v' \spr_{L^2(\Omega_-\cap S_n)} \\
&-\spl v, v' \spr_{L^2(\Omega)}
-\lambda \spl |\epsilon(\lambda)| v,v' \spr_{L^2(\Omega)}\\
&-2\lambda \epsilon_-|1-\omega^2_\epsilon/\lambda| \sum_{n=1}^{2N} \spl v,\chi_n R_n^{V,+} v' \spr_{L^2(\Omega_-\cap S_n)}\\
&-2\epsilon_-|1-\omega^2_\epsilon/\lambda| \sum_{n=1}^{2N} \spl \nabla w_0,(\nabla \chi_n) R_n^{W_0,+}w_0' \spr_{L^2(\Omega_-\cap S_n)}\\
&-\lambda\spl \epsilon(\lambda)\nabla w_0,T^{V,+}v' \spr_{L^2(\Omega)}
+\lambda\spl \epsilon(\lambda) v,\nabla T^{W_0,+}w_0' \spr_{L^2(\Omega)}
\end{align*}
for all $u,u'\in X$.
Hence $(T'(\lambda))^*A(\lambda)=A_1+A_2+A_3$.
Then we can follow the proof of \cite[Lemma~2]{BonnetBDCarvalhoCiarlet:18} line by line to obtain by means of a weighted Youngs inequality that $\spl A_2u,u \spr_X\geq c_2 \|\nabla w_0\|_{L^2(\Omega)}^2$ with a constant $c_2>0$ independent of $u\in X$.
The same technique allows to estimate $\spl A_1u,u \spr_X\geq c_1 \|v\|_X^2$ with a constant $c_1>0$ independent of $u\in X$.
Hence
\begin{align*}
\spl(A_1+A_2)u,u\spr_X\geq \min(c_1,c_2) (\|v\|_X^2+\|\nabla w_0\|_{L^2(\Omega)}^2).
\end{align*}
Since $V$ and $W$ are orthogonal and due to $\curl \nabla w_0=0$ it follows
\begin{align*}
\|v\|_X^2+\|\nabla w_0\|_{L^2(\Omega)}^2=\|u\|_X^2.
\end{align*}
Thus $A_1+A_2$ is coercive.
The operator $A_3$ is compact due to \eqref{eq:Vembeddings} and Assumption~\ref{ass:abstractReflection}.
\end{proof}

We continue with a discussion on the bijectivity of $T'(\lambda)$.
It can easily be seen that $(T^{V,\signv})^2=I_{X}$ and $(T^{W_0,\signw})^2=I_{W_0}$.
Since $\nabla T^{W_0,\signw} P_{W_0}$ maps into $W$, we further obtain that $(\nabla T^{W_0,\signw} P_{W_0})^2=I_W$.
However, $T^{V,\signv} P_V$ does not map into $V$ (for $\kappa_{\mu^{-1}}(\lambda)<0$).
For this reason we are not able to prove the bijectivity of $T'(\lambda)$ via the convenient way ($T'(\lambda)T'(\lambda)=I_X$).
However we can obtain a weaker result.

\begin{lemma}\label{lem:TprimeFred}
Let $\chi_n$ and $R_n^{V,\pm}, R_n^{W_{0,\pm}}$, $n=1,\dots,2N$ satisfy Assumptions~\ref{ass:abstractCutoff} and \ref{ass:abstractReflection} respectively.
Let $T'(\lambda)$ be defined as in \eqref{eq:DefTprime}.
Let $\lambda\in(0,\max(\omega^2_\mu,\omega^2_\epsilon))$ be such that $\kappa_{\mu^{-1}}(\lambda) \notin I_{R_V}$ and $\kappa_{\epsilon}(\lambda) \notin I_{R_{W_0}}$.
Then $T'(\lambda)$ is Fredholm with index zero.
Thus there exists a compact operator $K(\lambda)\in L(X)$ such that
\begin{align}\label{eq:DefT-crit}
T(\lambda):=T'(\lambda)+K(\lambda), \quad
\lambda\in(0,\max(\omega^2_\mu,\omega^2_\epsilon)) \colon
\kappa_{\mu^{-1}}(\lambda)\notin I_{R_V}, \kappa_{\epsilon}(\lambda)\notin I_{R_{W_0}}.
\end{align}
is bijective.
\end{lemma}
\begin{proof}
We abbreviate $T':=T'(\lambda)$ and compute
\begin{align*}
T'T'&=T^{V,\signv} P_V T^{V,\signv} P_V + P_{W}+ \nabla T^{W_0,\signw} P_{W_0} T^{V,\signv} P_V\\
&=T^{V,\signv} (I-P_W) T^{V,\signv} P_V + P_{W}+ \nabla T^{W_0,\signw} P_{W_0} T^{V,\signv} P_V\\
&=P_V -T^{V,\signv} P_W T^{V,\signv} P_V + P_{W}+ \nabla T^{W_0,\signw} P_{W_0} T^{V,\signv} P_V\\
&=I_X+(\nabla T^{W_0,\signw} -T^{V,\signv} \nabla)P_{W_0} T^{V,\signv} P_V.
\end{align*}
Thus if we can prove that $P_{W_0} T^{V,\signv} P_V$ is compact, it follows that $I_X-T'T'$ is compact too.
To see that the former operator is compact we recall that $P_{W_0}u$ solves
$\spl \nabla P_{W_0}u, \nabla w_0' \spr_{L^2(\Omega)}=\spl u, \nabla w_0' \spr_{L^2(\Omega)}$
for all $w_0'\in W_0=H^1_0(\Omega)$.
Hence $P_{W_0} T^{V,\signv} P_Vu$ solves
$\spl \nabla P_{W_0}T^{V,\signv} P_Vu, \nabla w_0' \spr_{L^2(\Omega)}=\spl T^{V,\signv} P_Vu, \nabla w_0' \spr_{L^2(\Omega)}$
for all $w_0'\in W_0=H^1_0(\Omega)$.
Let $E_{X,(L^2(\Omega))^2}\in L(X,(L^2(\Omega))^2)$ be the embedding operator.
Due to \eqref{eq:Vembeddings} and Assumption~\ref{ass:abstractReflection} it follows that $E_{X,(L^2(\Omega))^2} T^{V,\signv} P_V$ is compact.
Thus $P_{W_0} T^{V,\signv} P_V$ is compact and hence so is $K:=I_X-T'T'$.
So it follows with \cite[Lemma 2.5]{GohbergGoldbergKaashoek:03} that $T'$ is Fredholm and $2\ind T'=\ind(T'T')=\ind(I_X-K)=\ind I_X=0$.
The last claim follows with \cite[Corollary 2.4]{GohbergGoldbergKaashoek:03}.
\end{proof}

Let
\begin{align}\label{eq:DefLambdaR}
\begin{aligned}
\Lambda_{R_V,R_{W_0}}:=\,&\setC\setminus\setR_+\cup\big(\max(\omega_\mu^2,\omega_\epsilon^2),+\infty\big)\\
&\cup
\{\lambda\in (0,\max(\omega_\mu^2,\omega_\epsilon^2))\colon
\kappa_{\mu^{-1}}(\lambda) \notin I_{R_V} \text{ and } 
\kappa_{\epsilon}(\lambda) \notin I_{R_{W_0}}\}.
\end{aligned}
\end{align}

\begin{proposition}\label{prop:wTc-abstract}
Let $\chi_n$ and $R_n^{V,\pm}, R_n^{W_{0,\pm}}$, $n=1,\dots,2N$ satisfy Assumptions~\ref{ass:abstractCutoff} and \ref{ass:abstractReflection} respectively.
Let $A(\lambda)$ be defined as in \eqref{eq:DefA} and $T(\lambda)$ be defined as in \eqref{eq:DefT-identity}, \eqref{eq:DefT-classical} and \eqref{eq:DefT-crit}.
Then $A(\cdot)\colon \Lambda_{R_V,R_{W_0}}\to L(X)$ is weakly $T(\cdot)$-coercive and its resolvent set contains $\setC\setminus\setR_+$.
\end{proposition}
\begin{proof}
Follows from Section~\ref{sec:setting} and Lemmata~\ref{lem:TprimA-wc}, \ref{lem:TprimeFred}.
\end{proof}

\subsection{Explicit construction}\label{subsec:wTc-explicit}

In this subsection we construct cut-off functions $\chi_n$ and local reflection operators $R_n^{V,\pm}, R_n^{W_0,\pm}$ which satisfy Assumption~\ref{ass:abstractCutoff} and \ref{ass:abstractReflection} respectively.
To this end we follow \cite{BonnetBDCarvalhoCiarlet:18} very closely.
As therein we assume that all angles $\alpha_n\in2\pi\setQ$, $n=1,\dots,N$.
Let $R_n$, $R_n'$ be as specified in \cite[Appendix~A.1]{BonnetBDCarvalhoCiarlet:18} and $\chi_n$ be as specified in \cite[p.\ 816-817]{BonnetBDCarvalhoCiarlet:18}.
Then we set
\begin{align*}
R_n^{W_0,+}:=R_n \qquad\text{and}\qquad R_n^{W_0,-}:=R_n'.
\end{align*}
Thus $\chi_n$ satisfy Assumption~\ref{ass:abstractCutoff} and $R_n^{W_0,\pm}$ satisfy Assumption \ref{ass:abstractReflection}.
In addition we recall from \cite[Appendix~A.2]{BonnetBDCarvalhoCiarlet:18} that
\begin{align*}
\|R_n^{W_0,\pm}\|=\max \left( \frac{2\pi-\alpha_n}{\alpha_n},\frac{\alpha_n}{2\pi-\alpha_n} \right)
\end{align*}
for $n=1,\dots,N$ and $\|R_n^{W_0,\pm}\|=1$ for $n=N+1,\dots,2N$.
Thus
\begin{align}\label{eq:NormRW}
\|R^{W_0,\pm}\|=\max_{n=1,\dots,N} \max \left( \frac{2\pi-\alpha_n}{\alpha_n},\frac{\alpha_n}{2\pi-\alpha_n} \right).
\end{align}

\subsubsection{Details on the construction of $R_n^{W_0,\pm}$}
For the definition of $R_n^{V,\pm}$ we need to go into the details of the definitions of $R_n$ and $R_n'$.
For each corner $c_n$, $n=1,\dots,N$ let $p_{n,\pm}\in\setN$ be such that $\alpha_n=2\pi p_{n,+}/(p_{n,+}+p_{n,-})$ and $p_{n,+}+p_{n,-} \in 2\setN$.
(If $p_{n,+}+p_{n,-}$ is odd, replace the numbers by $2p_{n,+}$ and $2p_{n,-}$.)
We introduce the polar coordinates $(r_n,\theta_n)$ centered at $c_n$ such that $\Omega_+$ coincides locally with the cone $C_{\alpha_n}:=\{(r_n\cos\theta_n,r_n\sin\theta_n)\colon r_n>0,0<\theta_n<\alpha_n\}$.
Let $\beta_n:=2\pi/(p_{n,+}+p_{n,-})$ and $\calP\subset C_{\beta_n}$ be a bounded domain that coincides locally at the corner $c_n$ with $C_{\beta_n}$ and which is symmetric with respect to the angle $\beta/2$.
Denote the rotation $x=r_n\big(\cos\theta_n,\sin\theta_n\big)^\top\mapsto r_n\big(\cos(\theta_n+k\beta_n),\sin(\theta_n+k\beta_n)\big)^\top$ of angle $k\beta_n$ as $\calR_k$.
Then let $\Omega_{n,+}^k:=\calR_{k-1}\calP\subset\Omega_+$, $k=1,\dots,p_{n,+}$ and $\Omega_{n,-}^k:=\calR_{k-1+p_{n,+}}\calP\subset\Omega_-$, $k=1,\dots,p_{n,-}$.
Subsequently define $\Omega_n$ as the interior of the closure of $\bigcup_{k=1}^{p_{n,+}} \Omega_{n,+}^k \cup \bigcup_{k=1}^{p_{n,-}} \Omega_{n,-}^k$.
Then $R_n^{W_0,\pm}$ is defined piece-wise
\begin{align*}
(R_n^{W_0,\pm}w_0)|_{\Omega_{n,\mp}^k}=\sum_{m=1}^{p_{n,\pm}} G^{W_0}_{n,\pm,k,m} w_0,
\quad k=1,\dots,p_{n,\mp}
\end{align*}
with operators
\begin{align*}
G^{W_0}_{n,\pm,k,m}w_0=w_0\circ\phi_{n,\pm,k,m}, \quad k=1,\dots,p_{n,\mp}, \quad m=1,\dots,p_{n,\pm}
\end{align*}
which are itself defined in terms of affine transformations
\begin{align*}
\phi_{n,\pm,k,m}(x)=F_{n,\pm,k,m} x + y_{n,\pm,k,m}
\end{align*}
with orthonormal matrices $F_{n,\pm,k,m}\in\setR^{2\times 2}$ and $y_{n,\pm,k,m} \in \setR^2$.
In particular $\phi_{n,\pm,k,m}$ maps $\Omega_{n,+}^m$ onto $\Omega_{n,-}^k$.
For $n=N+1,\dots,2N$ one introduces a neighborhood $\Omega_n$ of $e_n$ with $\ol{S_n}\subset\Omega_n$ which is symmetric with respect to $e_n$ and $\Omega_{n,\pm}:=\Omega_n\cap\Omega_\pm$.
The local reflection operators are defined with a single transformation
\begin{align*}
(R_n^{W_0,\pm}w_0)|_{\Omega_{n,\mp}}= G^{W_0}_{n,\pm} w_0
\end{align*}
whereby again
\begin{align*}
G^{W_0}_{n,\pm}w_0=w_0\circ\phi_{n,\pm}
\end{align*}
with affine transformations
\begin{align*}
\phi_{n,\pm}(x)=F_{n,\pm} x + y_{n,\pm},
\end{align*}
with orthonormal matrices $F_{n,\pm}\in\setR^{2\times 2}$ and $y_{n,\pm} \in \setR^2$,
whereby $\phi_{n,\pm}$ maps $\Omega_{n,\pm}$ onto $\Omega_{n,\mp}$.
At last each local reflection $(R_n^{W_0,\pm}w_0)|_{\Omega_\mp\cap\Omega_n}$, $n=1,\dots,2N$ is extended to $\Omega_\mp$ in an arbitrary way (only products $\chi_n R_n^{W_0,\pm}$ appear in our analysis).
The analysis in \cite{BonnetBDCarvalhoCiarlet:18} shows that for $n=1,\dots,N$ one can estimate
\begin{align*}
\|\chi_n^{1/2}\nabla R_n^{W_0,\pm} w_0\|_{L^2(\Omega_\mp\cap S_n)}
\leq | M_{n,\pm} W |_2
&\leq | M_{n,\pm}^\top M_{n,\pm} |_2 |W|_2\\
&= | M_{n,\pm}^\top M_{n,\pm} |_2 \|\chi_n^{1/2}\nabla w_0\|_{L^2(\Omega_\pm\cap S_n)}
\end{align*}
with matrices $M_{n,\pm}\in\setR^{p_{n,-},p_{n,+}}$, vectors
\begin{align*}
W:=(\|\chi_n^{1/2}\nabla w_0\|_{L^2(\Omega_{n,+}^1)},\dots,\|\chi_n^{1/2}\nabla w_0\|_{L^2(\Omega_{n,+}^{p_{n,+}})})^\top
\end{align*}
and the euclidian norm $|\cdot|_2$.
The core of the analysis in \cite{BonnetBDCarvalhoCiarlet:18} is the estimate $| M_{n,\pm}^\top M_{n,\pm} |_2$ $\leq$ $\max \left( \frac{2\pi-\alpha_n}{\alpha_n},\frac{\alpha_n}{2\pi-\alpha_n} \right)$ and thus $\|R_n^{W_0,\pm}\| \leq \max \left( \frac{2\pi-\alpha_n}{\alpha_n},\frac{\alpha_n}{2\pi-\alpha_n} \right)$.
For $n=N+1,\dots,2N$ the estimate $\|R_n^{W_0,\pm}\| \leq 1$ follows with less effort.
Subsequently one can choose $w_0$ as in \cite{BonnetBDCarvalhoCiarlet:18} to find that indeed equality holds in the norm estimates.
\sloppy

\subsubsection{Construction of $R_n^{V,\pm}$}
To define $R_n^{V,\pm}$ we recall the transformation rule (see e.g.\ \cite[Lemma 4.15]{Zaglmayr:06})
\begin{align*}
J (\curl u)\circ\phi = \curl\big( F^\top u\circ\phi \big),
\end{align*}
whereby $\phi$ is a domain transformation, $F:=D\phi$ and $J:=\det F$.
Then we set
\begin{align*}
G^V_{n,\pm,k,m}u:=F_{n,\pm,k,m}^\top u\circ\phi_{n,\pm,k,m}, \quad k=1,\dots,p_{n,\mp}, \quad m=1,\dots,p_{n,\pm}
\end{align*}
and
\begin{align*}
(R_n^{V,\pm}u)|_{\Omega_{n,\mp}^k}:=\sum_{m=1}^{p_{n,\pm}} G^V_{n,\pm,k,m} u, \quad k=1,\dots,p_{n,\mp}
\end{align*}
for $n=1,\dots,N$ and
\begin{align*}
G^V_{n,\pm}u:=F_{n,\pm}^\top u\circ\phi_{n,\pm}
\end{align*}
and
\begin{align*}
(R_n^{V,\pm}u)|_{\Omega_{n,\mp}}:=G^V_{n,\pm} u
\end{align*}
for $n=N+1,\dots,2N$ in analogy to $R_n^{W_0,\pm}$.
It follows that $R_n^{V,\pm}$ satisfy the matching condition on the traces in Assumption~\ref{ass:abstractReflection}.
To see that $R_n^{V,\pm} P_V \in L(X_\pm,H^t(\Omega_\mp\cap\Omega_n)$ for a $t\in(0,1/2)$ we recall \eqref{eq:Vembeddings} that $V$ continuously embeds into $(H^{1/2}(\Omega))^2$.
From the construction of $R_n^{V,\pm}$ it follows that $R_n^{V,\pm}P_Vu$ is in $(H^{1/2}(\Omega_{n,\mp}^k))^2$ for each $k=1,\dots,p_{,-}$ and thus in $(H^{t}(\Omega_{n,\mp}^k))^2$ for each $k=1,\dots,p_{,-}$ and $t\in(0,1/2)$.
Since for $t\in(0,1/2)$ the piece-wise regularity already implies the global regularity,
$R_n^{V,\pm}$ satisfy Assumption~\ref{ass:abstractReflection}.
From the definition of $R_n^{V,\pm}$ it follows that the analysis of \cite[Appendix~A.2]{BonnetBDCarvalhoCiarlet:18} needs only to be repeated to compute the norms $\|R_n^{V,\pm}\|$, e.g.\ for $n=1,\dots,N$ one can estimate
\begin{align*}
\|\chi_n^{1/2}\curl R_n^{V,\pm} u\|_{L^2(\Omega_\mp\cap S_n)}
\leq | M_{n,\pm} U |_2
&\leq | M_{n,\pm}^\top M_{n,\pm} |_2 |U|_2\\
&= | M_{n,\pm}^\top M_{n,\pm} |_2 \|\chi_n^{1/2}\curl u\|_{L^2(\Omega_\pm\cap S_n)}
\end{align*}
with matrices $M_{n,\pm}\in\setR^{p_{n,-},p_{n,+}}$ as hitherto and
\begin{align*}
U:=(\|\chi_n^{1/2}\curl u\|_{L^2(\Omega_{n,+}^1)},\dots,\|\chi_n^{1/2}\curl u\|_{L^2(\Omega_{n,+}^{p_{n,+}})})^\top.
\end{align*}
Thus we obtain
\begin{align*}
\|R_n^{V,\pm}\|=\max \left( \frac{2\pi-\alpha_n}{\alpha_n},\frac{\alpha_n}{2\pi-\alpha_n} \right)
\end{align*}
for $n=1,\dots,N$ and $\|R_n^{V,\pm}\|=1$ for $n=N+1,\dots,2N$.
Hence
\begin{align}\label{eq:NormRV}
\|R^{V,\pm}\|=\max_{n=1,\dots,N} \max \left( \frac{2\pi-\alpha_n}{\alpha_n},\frac{\alpha_n}{2\pi-\alpha_n} \right).
\end{align}
Thus we define the critical interval
\begin{align*}
I_c:=[-I_\alpha,-1/I_\alpha], \qquad I_\alpha:=\max_{n=1,\dots,N} \max \left( \frac{2\pi-\alpha_n}{\alpha_n},\frac{\alpha_n}{2\pi-\alpha_n} \right)
\end{align*}
and
\begin{align}\label{eq:DefLambdaAlpha}
\Lambda_{\alpha}:=\{\lambda\in\setC\colon \lambda\neq 0,\omega_\mu^2,\omega_\epsilon^2
\text{ and }
\kappa_{\mu^{-1}}(\lambda), \kappa_{\epsilon}(\lambda) \notin I_c\}.
\end{align}

\subsection{Summary of results}\label{subsec:wTc-results}

We combine the results from Subsections~\ref{subsec:wTc-abstract} and \ref{subsec:wTc-explicit}.

\begin{theorem}\label{thm:wTc-explicit}
Assume that all the corners' angles of the interface $\Sigma$ belong to $2\pi\setQ$.
Let $A(\cdot)$ be as defined in \eqref{eq:DefA}, $\chi_n$, $R_n^{V,\pm}$, $R_n^{W_0,\pm}$ be as defined in Subsection~\ref{subsec:wTc-explicit} and $T(\cdot)$ be as defined in \eqref{eq:DefT-identity}, \eqref{eq:DefT-classical} and \eqref{eq:DefT-crit}.
Then $A(\cdot)\colon \Lambda_{\alpha} \to L(X)$ is weakly $T(\cdot)$-coercive and its resolvent set contains $\setC\setminus\setR_+$.
\end{theorem}
\begin{proof}
Follows from Proposition~\ref{prop:wTc-abstract} and Subsection~\ref{subsec:wTc-explicit}.
\end{proof}

\begin{corollary}\label{cor:SpecProp}
Let the assumptions of Theorem~\ref{thm:wTc-explicit} be satisfied.
Then the spectrum of $A(\cdot)\colon \Lambda_{\alpha} \to L(X)$ is discrete, has no accumulation points in $\Lambda_{\alpha}$ and every point in the spectrum is an eigenvalue with finite algebraic multiplicity.
\end{corollary}
\begin{proof}
Follows from Theorem~\ref{thm:wTc-explicit} and the analytic Fredholm theorem.
\end{proof}

\begin{remark}
The studies in the preceding articles \cite{BonnetBDChesnelCiarlet:12,BonnetBDChesnelCiarlet:14b,BonnetBDChesnelCiarlet:14a,BonnetBDCarvalhoCiarlet:18} considered the case of dispersionless $\mu$, $\epsilon$ and were not able to obtain that the resolvent set of the corresponding operator function is non-empty.
Though for dispersive $\mu$, $\epsilon$ it is no problem to obtain coercivity for $\lambda\in\setC\setminus\setR_+$.
\end{remark}

\section{T-compatible approximation}\label{sec:Tcomp}

In this section we analyze the approximation of $A(\cdot)$ by Galerkin methods.
To this end we employ the framework of $T(\cdot)$-compatible approximations \cite{Halla:19Tcomp}.
In Subsection~\ref{subsec:Tcomp-abstract} we assume that Assumptions~\ref{ass:abstractCutoff}, \ref{ass:abstractReflection} are satisfied and formulate four additional Assumptions \ref{ass:proj}, \ref{ass:localpol}, \ref{ass:localconst}, \ref{ass:Rconform} on the Galerkin spaces.
Under this pretext we prove for such approximations their $T(\cdot)$-compatibility.
In Subsection~\ref{subsec:Tcomp-explicit} we discuss that h-finite element methods with N\'ed\'elec-elements and locally $R$-conform meshes satisfy the former assumptions.
In Subsection~\ref{subsec:Tcomp-results} we conclude and formulate our results.

\subsection{Abstract framework}\label{subsec:Tcomp-abstract}

First we formulate our assumptions.

\begin{assumption}[$L^2$-uniformly bounded commuting cochain projections]\label{ass:proj}
For the de Rham Complex
\begin{align*}
W_0=H^1_0(\Omega)
\xrightarrow{\nabla} X=H_0(\curl;\Omega)
\xrightarrow{\curl} J:=L^2(\Omega).
\end{align*}
let for each $n\in\setN$ 
\begin{align*}
W_{0,n}
\xrightarrow{\nabla} X_n
\xrightarrow{\curl} J_n
\end{align*}
be a subcomplex with finite dimensional spaces $W_{0,n}$, $X_n$, $J_n$ such that the corresponding orthogonal projections $P_{W_{0,n}}$, $P_{X_n}$, $P_{J_n}$ converge point-wise to the identity for $n\to\infty$.
We assume the existence of $L^2$-uniformly bounded commuting projections $\poop_{W_{0,n}}$, $\pi_{X_n}$, $\pi_{J_n}$ onto the subcomplex.
By that we mean that the following commuting diagram
\begin{align*}
\begin{array}{ccccc}
W_0 
&\xrightarrow{\nabla} & X 
&\xrightarrow{\curl} & J \\
\poop_{W_{0,n}}\downarrow && \poop_{X_n}\downarrow && \poop_{J_n}\downarrow \\
W_{0,n}
&\xrightarrow{\nabla}& X_n
&\xrightarrow{\curl} & J_n
\end{array}
\end{align*}
is satisfied, the projections are surjective and of the form
\begin{align*}
\pi_{W_{0,n}} = \tilde\pi_{W_{0,n}} E_{W_0}, \qquad
\pi_{X_n} = \tilde\pi_{X_n} E_X, \qquad
\pi_{J_n} = \tilde\pi_{J_n} E_J,
\end{align*}
with embedding operators $E_{W_0}\in L\big(W_0,L^2(\Omega)\big)$, $E_X\in L\big(X,(L^2(\Omega))^2\big)$ and $E_J\in L\big(J,L^2(\Omega)\big)$
and projection operators $\tilde\pi_{W_{0,n}} \in L\big(L^2(\Omega),W_{0,n}\big)$, $\tilde\pi_{X_n}\in L\big((L^2(\Omega))^2,X_n\big)$, $\tilde\pi_{J_n}\in L\big(L^2(\Omega),J_n\big)$ such that
\begin{align*}
\sup_{n\in\setN} \big\{
\|\tilde\pi_{W_{0,n}}\|_{L(L^2(\Omega))},\,
\|\tilde\pi_{X_n}\|_{L((L^2(\Omega))^2)},\,
\|\tilde\pi_{J_n}\|_{L(L^2(\Omega))}
\big\} < \infty.
\end{align*}
\end{assumption}

\begin{assumption}[local projection operators]\label{ass:localpol}
Let Assumption~\ref{ass:proj} be satisfied and
$X_{\pm,n}:=\{u|_{\Omega_\pm}\colon u\in X_n\}\subset X_\pm$,
$W_{0,\pm,n}:=\{w_0|_{\Omega_\pm}\colon w_0\in W_{0,n}\}\subset W_{0,\pm,}$.
Let $\big(h_n \big)_{n\in\setN}\in (\setR^+)^\setN$ be a non-increasing sequence with $\lim_{n\in\setN} h_n =0$ and $\CC, \CR >1$.
There exist projection operators $\ipop_{X_{\pm,n}}\in L(X_\pm,X_{\pm,n})$, $\ipop_{W_{0,\pm,n}}\in L(W_{0,\pm},W_{0,\pm,n})$, $n\in\setN$, that act locally in the following sense:
for $n\in\setN$, $s \in \{1,2\}$, $x_0\in\Omega$, if $u\in X$, $w_0\in W_0$ and $u|_{\Omega_\pm\cap B_{\CR  h_n }(x_0)} \in H^{s-1}(\curl;\Omega_\pm\cap B_{\CR  h_n }(x_0))$, $w_0|_{B_{\CR  h_n }(x_0)} \in H^s(B_{\CR  h_n }(x_0))$, then
\begin{align*}
\|u-\ipop_{X_{\pm,n}} u\|_{H(\curl;(\Omega_\pm\cap B_{h_n }(x_0))} &\leq \CC h_n ^{s-1} \|u\|_{H^{s-1}(\curl;\Omega_\pm\cap B_{\CR  h_n }(x_0))},\\
\|w_0-\ipop_{W_{0,\pm,n}} w_0\|_{H^1(\Omega_\pm\cap B_{h_n }(x_0))} &\leq \CC h_n ^{s-1} \|w_0\|_{H^s(\Omega_\pm\cap B_{\CR  h_n }(x_0))}.
\end{align*}
In addition, if $\tau\cdot u|_\Sigma \in \{\tau\cdot u'|_\Sigma\colon u'\in X_n\}$,
$w_0|_\Sigma \in \{w_0'|_\Sigma\colon w_0'\in W_{0,n}\}$, then
$\tau\cdot \ipop_{X_{\pm,n}}u|_\Sigma = \tau\cdot u|_\Sigma$,
$\ipop_{W_{0,\pm,n}}w_0|_\Sigma = w_0|_\Sigma$.
\end{assumption}

\begin{assumption}[local constants]\label{ass:localconst}
Let Assumption~\ref{ass:proj} be satisfied.
For each $D\subset\Omega$ which is compact in $\Omega$ exists $n_0>0$ such that for each $n\in\setN, n>n_0$ there exists $u_{D,n,1}, u_{D,n,2}\in X_n$, $w_{0,D,n}\in W_{0,n}$ with $u_{D,n,1}|_D=(1,0)^\top$, $u_{D,n,2}|_D=(0,1)^\top$, $w_{0,D,n}|_D=1$.
\end{assumption}

\begin{assumption}[local $R$-conformity]\label{ass:Rconform}
Let Assumptions~\ref{ass:abstractCutoff}, \ref{ass:abstractReflection} and \ref{ass:proj} be satisfied.
For all $n\in\setN$, $m=1,\dots,2N$ and $u\in X_n$, $w_0\in W_{0,n}$ it holds
\begin{align*}
(R_m^{V,\pm}u)|_{\Omega_{n,\mp}} &\in \{u|_{\Omega_\mp\cap\Omega_m}\colon u\in X_n\},\\
(R_m^{W_0,\pm}w_0)|_{\Omega_{n,\mp}} &\in \{w_0|_{\Omega_\mp\cap\Omega_m}\colon w_0\in W_{0,n}\}.
\end{align*}
\end{assumption}
For approximations which satisfy Assumptions~\ref{ass:proj}, \ref{ass:localpol}, \ref{ass:localconst} and \ref{ass:Rconform} we define
\begin{align*}
T_n^{V,\pm}u&:=\left\{\begin{array}{ll}
\phantom{-}u|_{\Omega_+}-(1\mp1)\ipop_{X_{+,n}}\sum_{n=1}^{2N}\chi_n R_n^{V,-}u|_{\Omega_-}, &\text{in }\Omega_+,\\
-u|_{\Omega_-}+(1\pm1)\ipop_{X_{-,n}}\sum_{n=1}^{2N}\chi_n R_n^{V,+}u|_{\Omega_+}, &\text{in }\Omega_-,
\end{array}\right.,\\
T_n^{W_0,\pm}w_0&:=\left\{\begin{array}{ll}
\phantom{-}w_0|_{\Omega_+}-(1\mp1)\ipop_{W_{0,+,n}}\sum_{n=1}^{2N}\chi_n R_n^{W_0,-}w_0|_{\Omega_-}, &\text{in }\Omega_+,\\
-w_0|_{\Omega_-}+(1\pm1)\ipop_{W_{0,-,n}}\sum_{n=1}^{2N}\chi_n R_n^{W_0,+}w_0|_{\Omega_+}, &\text{in }\Omega_-,
\end{array}\right..
\end{align*}
In addition we set $T_n^{V,0}:=I_X$ and $T_n^{W_0,0}:=I_{W_0}$.
Then we consider the composed operator
\begin{subequations}
\begin{align}\label{eq:DefTnprime}
T_n'(\lambda):=T_n^{V,\signv} \poop_{X_n} P_V - \nabla T_n^{W_0,\signw} \poop_{W_{0,n}} P_{W_0}
\end{align}
with
\begin{align}
\begin{aligned}
\begin{array}{ll}
\signv=0, &\text{if } \kappa_{\mu^{-1}}(\lambda)>0,\\
\signv=-, &\text{if } \kappa_{\mu^{-1}}(\lambda)<-1,\\
\signv=+, &\text{if } \kappa_{\mu^{-1}}(\lambda)\in(-1,0),
\end{array}
\text{and}\quad
\begin{array}{ll}
\signw=0, &\text{if } \kappa_\epsilon(\lambda)>0,\\
\signw=-, &\text{if } \kappa_\epsilon(\lambda)<-1,\\
\signw=+, &\text{if } \kappa_\epsilon(\lambda)\in(-1,0),
\end{array}.
\end{aligned}
\end{align}
\end{subequations}
Further let
\begin{align}\label{eq:DefTn}
T_n(\lambda):=\left\{\begin{array}{ll}
I_X, &\lambda\in\setC\setminus\setR_+,\\
\poop_{X_n} T(\lambda)=\poop_{X_n} P_V - \poop_{X_n} P_W, &\lambda>\max(\omega^2_\mu,\omega^2_\epsilon),\\
T_n'(\lambda)+P_{X_n} K(\lambda), &\lambda\in(0,\max(\omega^2_\mu,\omega^2_\epsilon)),
\end{array}\right..
\end{align}
It follows from Assumptions~\ref{ass:proj}, \ref{ass:localpol} that indeed $T_n(\lambda)\in L(X_n)$.
Now that we have constructed $T_n(\lambda)$ we want to prove $\lim_{n\to\infty} \|T(\lambda)-T_n(\lambda)\|_{L(X_n,X)}=0$.
We split the proof into several lemmata and summarize in Theorem~\ref{thm:Tcomp}.
The following result is quite standard.

\begin{lemma}\label{lem:piPV}
Let Assumption~\ref{ass:proj} be satisfied.
Then there hold
\begin{align*}
\lim_{n\to\infty} \| (I_X-\poop_{X_n}) P_V \|_{L(X_n,X)}=
\lim_{n\to\infty} \| (I_X-\poop_{X_n}) P_W \|_{L(X_n,X)}=0.
\end{align*}
\end{lemma}
\begin{proof}
See \cite{ArnoldFalkWinther:10}.
We state a proof for completeness.
It follows from $\curl\poop_{X_n}P_W=\poop_{J_n}\curl P_W=0$ that $\| (I_X-\poop_{X_n}) P_W \|_{L(X_n,X)}=\| (I_X-\poop_{X_n}) P_W \|_{L(X_n,(L^2(\Omega))^2)}$.
Since $\poop_{X_n}$ is a projection onto $X_n$ and $P_W=I_X-P_V$ it follows $(I_X-\poop_{X_n}) P_W |_{X_n}$ $=$ $-$ $(I_X-\poop_{X_n}) P_V |_{X_n}$.
Thus
\begin{align*}
\|(I_X-\poop_{X_n}) P_V \|_{L(X_n,X)}&=\| (I_X-\poop_{X_n}) P_W \|_{L(X_n,X)}\\
&= \| (I_X-\poop_{X_n}) P_W \|_{L(X_n,(L^2(\Omega))^2)}\\
&= \|(I_X-\poop_{X_n}) P_V \|_{L(X_n,(L^2(\Omega))^2)}\\
&\leq \| (I_{(L^2(\Omega))^2}-\tilde\poop_{X_n}) E_X P_V \|_{L(X,(L^2(\Omega))^2)}.
\end{align*}
Since $E_X P_V$ is compact and $(I_{(L^2(\Omega))^2}-\tilde\poop_{X_n})$ converges point-wise to zero the claim is proven.
\end{proof}

\begin{lemma}\label{lem:Vpoop}
Let Assumptions~\ref{ass:abstractCutoff}, \ref{ass:abstractReflection} and \ref{ass:proj}, \ref{ass:localpol}, \ref{ass:localconst}, \ref{ass:Rconform} be satisfied.
Then for each $\signv\in\{0,-,+\}$ it holds
$\lim_{n\to\infty} \| T^{V,\signv} (I_X-\poop_{X_n}) P_V \|_n=0$.
\end{lemma}
\begin{proof}
We estimate
\begin{align*}
\| T^{V,\signv} (I_X-\poop_{X_n}) P_V \|_{L(X_n,X)}
\leq \|T^{V,\signv}\|_{L(X)} \| (I_X-\poop_{X_n}) P_V \|_{L(X_n,X)}
\end{align*}
and hence the claim follows with Lemma \ref{lem:piPV}.
\end{proof}

\begin{lemma}\label{lem:Wpoop}
Let Assumptions~\ref{ass:abstractCutoff}, \ref{ass:abstractReflection} and \ref{ass:proj}, \ref{ass:localpol}, \ref{ass:localconst}, \ref{ass:Rconform} be satisfied.
Then for each $\signw\in\{0,-,+\}$ it holds
$\lim_{n\to\infty}\| T^{W_0,\signw} (I_{W_0}-\poop_{W_{0,n}}) P_{W_0} \|_{L(X_n,W_0)}$ $=$ $0$.
\end{lemma}
\begin{proof}
We estimate
\begin{align*}
\| T^{W_0,\signw} (I_{W_0}-\poop_{W_{0,n}}) & P_{W_0} \|_{L(X_n,W_0)}\\
&\leq \|T^{W_0,\signw}\|_{L(W_0)} \|(I_{W_0}-\poop_{W_{0,n}}) P_{W_0} \|_{L(X_n,W_0)}.\end{align*}
Since on $W_0$ the norms $\|\cdot\|_{H^1(\Omega)}$ and$\|\nabla \cdot\|_{L^2(\Omega)}$ are equivalent, we estimate further
\begin{align*}
\|(I_{W_0}-\poop_{W_{0,n}}) P_{W_0} \|_{L(X_n,W_0)}
\lesssim \|\nabla (I_{W_0}-\poop_{W_{0,n}}) P_{W_0} \|_{L(X_n,X)}.
\end{align*}
Now we reformulate $\nabla (I_{W_0}-\poop_{W_{0,n}}) P_{W_0}|_{X_n}=(I_X-\poop_{X_n})P_W|_{X_n}$ and apply Lemma \ref{lem:piPV}.
\end{proof}

For the proofs of Lemmata \ref{lem:Vipop} and \ref{lem:Wipop} we adapt the discrete commutator technique \cite{Bertoluzza:99}.

\begin{lemma}\label{lem:Vipop}
Let Assumptions~\ref{ass:abstractCutoff}, \ref{ass:abstractReflection} and \ref{ass:proj}, \ref{ass:localpol}, \ref{ass:localconst}, \ref{ass:Rconform} be satisfied.
Then for each $\signv\in\{0,-,+\}$ it holds
$\lim_{n\to\infty} \| (T^{V,\signv}-T_n^{V,\signv}) \poop_{X_n} P_V \|_{L(X_n,X)}=0$.
\end{lemma}
\begin{proof}
For $\signv=0$ it holds $(T^{V,\signv}-T_n^{V,\signv}) \poop_{X_n} P_V=0$.
Thus we consider the remaining cases $\signv=\pm$.
We slightly adapt the technique introduced in \cite{Bertoluzza:99}.
We compute
\begin{align*}
\| (T^{V,\pm}-T_n^{V,\pm}) &\poop_{X_n} P_V \|_{L(X_n,X)}\\
&=2\sup_{u\in X_n, \|u\|_X=1} \| \sum_{m=1}^{2N} (1-\ipop_{X_{\mp,n}})\chi_m R_m^{V,\pm}\poop_{X_n}P_Vu\|_{H(\curl;\Omega_\mp)} \\
&\leq 2\sum_{m=1}^{2N} \sup_{u\in X_n, \|u\|_X=1} \| (1-\ipop_{X_{\mp,n}})\chi_m R_m^{V,\pm}\poop_{X_n}P_Vu\|_{H(\curl;\Omega_\mp)}
\end{align*}
Now for each $m=1,\dots,2N$ we proceed as follows.
From Assumptions~\ref{ass:abstractCutoff}, \ref{ass:abstractReflection} and \ref{ass:localpol} it follows as in the proof of \cite[Theorem 5.5]{Halla:20PML} that there exist $n_0>0$ and a Lipschitz domain $B\subset\Omega_\mp\cap\Omega_m$ with $S_m\subset B$ and $\ol{B}\subset\Omega_m$ such that
\begin{align*}
&\sup_{u\in X_n, \|u\|_X=1} \| (1-\ipop_{X_{\mp,n}})\chi_m R_m^{V,\pm}\poop_{X_n}P_Vu\|_{H(\curl;\Omega_\mp)}\\
=&\sup_{u\in X_n, \|u\|_X=1} \| (1-\ipop_{X_{\mp,n}})\chi_m R_m^{V,\pm}\poop_{X_n}P_Vu\|_{H(\curl;B)}
\end{align*}
for all $n>n_0$.
For each $n>n_0$ we consider a collection of balls $\{B_{h_n }(x)\colon x\in Z\}$ such that $B\subset \bigcup_{x\in Z} B_{h_n }(x)$ and such that any point $y\in \Omega$ belongs to at most $M\in\setN$ (with $M$ independent of $n\in\setN$) balls of the collection $\{B_{\CR  h_n }(x)\colon x\in Z\}$ and such that there exists $n_1>n_0$ with
$\tilde B:=\bigcup_{x\in Z} B_{\CR  h_n }(x))\cap\Omega_\mp \subset \Omega_m\cap\Omega_\mp$ for all $n>n_1$.
It follows from the construction of the covering that there exists a constant $C_1>0$ such that
\begin{align*}
\sum_{x\in Z} \|u\|^2_{H^s(\curl;B_{\CR  h_n }(x))} \leq C_1 \|u\|^2_{H^s(\curl;\tilde B)},\quad
\sum_{x\in Z} \|v\|^2_{H^s(B_{\CR  h_n }(x))} \leq C_1 \|u\|^2_{H^s(\tilde B)}
\end{align*}
for $s\in[0,2]$ and all $u\in H^s(\curl;\tilde B)$, $v\in H^s(\tilde B)$ and $n>n_1$.
We abbreviate $v:=R_m^{V,\pm} P_V u$ and choose $v_x\in X_{\mp,n}$ for each $x\in Z$ such that $v_x|_{B_{\CR  h_n }(x)\cap \Omega_\mp}$ is constant and
\begin{align*}
\|v_x\|_{L^2(B_{\CR  h_n }(x)\cap \Omega_\mp)} &\leq \|v\|_{L^2(B_{\CR  h_n }(x)\cap \Omega_\mp)},\\
\|v-v_x\|_{L^2(B_{\CR  h_n }(x)\cap \Omega_\mp)} &\leq C_2 h_n^t \|v\|_{H^t(B_{\CR  h_n }(x)\cap \Omega_\mp)},
\end{align*}
with a constant $C_2>0$.
Thereby $t\in(0,1/2)$ is as in Assumption~\ref{ass:abstractReflection}.
This is possible due to the Bramble-Hilbert Lemma and interpolation techniques.
We estimate
\begin{subequations}
\begin{align}
\nonumber
\| (1-&\ipop_{X_{\mp,n}})\chi_m R_m^{V,\pm}\poop_{X_n}P_Vu\|_{H(\curl;B)}\\
\nonumber
&\leq \sum_{x\in Z} \| (1-\ipop_{X_n}) \chi_m R_m^{V,\pm} \poop_{X_n} P_V u\|_{H(\curl;B_{h_n }(x)\cap\Omega_\mp)}\\
\nonumber
&=\sum_{x\in Z} \| (1-\ipop_{X_{\mp,n}}) (\chi_m-\chi_m(x)) R_m^{V,\pm} \poop_{X_n} P_V u\|_{H(\curl;B_{h_n }(x)\cap\Omega_\mp)}\\
\nonumber
&\leq \sum_{x\in Z} \| (1-\ipop_{X_{\mp,n}}) (\chi_m-\chi_m(x)) R_m^{V,\pm} P_V u\|_{H(\curl;B_{h_n }(x)\cap\Omega_\mp)}\\
\nonumber
&+\sum_{x\in Z} \| (1-\ipop_{X_{\mp,n}}) (\chi_m-\chi_m(x)) R_m^{V,\pm} (I_X-\poop_{X_n}) P_V u\|_{H(\curl;B_{h_n }(x)\cap\Omega_\mp)} \\
\nonumber
&\leq \sum_{x\in Z} \| (1-\ipop_{X_{\mp,n}}) (\chi_m-\chi_m(x)) (v-v_x)\|_{H(\curl;B_{h_n }(x)\cap\Omega_\mp)} \\
\nonumber
&+\sum_{x\in Z} \| (1-\ipop_{X_{\mp,n}}) (\chi_m-\chi_m(x)) v_x\|_{H(\curl;B_{h_n }(x)\cap\Omega_\mp)} \\
\nonumber
&+\sum_{x\in Z} \| (1-\ipop_{X_{\mp,n}}) (\chi_m-\chi_m(x)) R_m^{V,\pm} (I_X-\poop_{X_n}) P_V u\|_{H(\curl;B_{h_n }(x)\cap\Omega_\mp)} \\
\label{eq:termi}
&\leq \CC \sum_{x\in Z} \| (\chi_m-\chi_m(x)) (v-v_x)\|_{H(\curl;B_{\CR  h_n }(x)\cap\Omega_\mp)} \\
\label{eq:termii}
&+\CC h_n  \sum_{x\in Z} \| (\chi_m-\chi_m(x)) v_x\|_{H^1(\curl;B_{\CR  h_n }(x)\cap\Omega_\mp)} \\
\label{eq:termiii}
&+\CC \sum_{x\in Z} \| (\chi_m-\chi_m(x)) R_m^{V,\pm} (I_X-\poop_{X_n}) P_V u\|_{H(\curl;\CR  B_{h_n }(x)\cap\Omega_\mp)}.
\end{align}
\end{subequations}
The above equality holds, because $\big( \ipop_{X_{\mp,n}} \chi_m(x) R_m^{V,\pm} \poop_{X_n} P_V u\big)|_{B_{h_n }(x)\cap\Omega_\mp}$ depends only on
\begin{align*}
\big(\chi_m(x) R_m^{V,\pm} \poop_{X_n} P_V u\big)|_{B_{\CR h_n }(x)\cap\Omega_\mp} \in \{u|_{B_{\CR h_n }(x)\cap\Omega_\mp}\colon u\in X_n\}
\end{align*}
and because $\ipop_{X_{\mp,n}}$ is a projection.
We continue to estimate
\begin{align*}
\eqref{eq:termii} &\lesssim 
\CC h_n  \sum_{x\in Z} \| (\chi_m-\chi_m(x)) v_x\|_{H^2(B_{\CR  h_n }(x)\cap\Omega_\mp)} \\
&\leq \CC h_n  \sum_{x\in Z} \| \chi_m-\chi_m(x)\|_{W^{2,\infty}(B_{\CR  h_n }(x)\cap\Omega_\mp)} \| v_x\|_{L^2(B_{\CR  h_n }(x)\cap\Omega_\mp)} \\
&\leq \CC h_n  2 \| \chi_m\|_{W^{2,\infty}(\Omega)} \sum_{x\in Z}  \| v_x\|_{L^2(B_{\CR  h_n }(x)\cap\Omega_\mp\cap\Omega_\mp)} \\
&\leq \CC h_n  2 \| \chi_m\|_{W^{2,\infty}(\Omega)} \sum_{x\in Z}  \| v\|_{L^2(B_{\CR  h_n }(x))} \\
&\leq \CC h_n  2 \| \chi_m\|_{W^{2,\infty}(\Omega)} C_1 \|v\|_{L^2(\tilde B)} \\
&\leq \CC h_n  2 \| \chi_m\|_{W^{2,\infty}(\Omega)} C_1 \|R_m^{V,\pm}\|_{L(X_\pm,X_\mp)} \|u\|_X.
\end{align*}
Further
\begin{align*}
\eqref{eq:termiii} &\lesssim
\CC \sum_{x\in Z} \|\chi_m-\chi_m(x)\|_{W^{1,\infty}(B_{\CR  h_n }(x)\cap \Omega_\mp)} \|R_m^{V,\pm} (I_X-\poop_{X_n}) P_V u\|_{H(\curl; B_{\CR  h_n }(x)\cap \Omega_\mp)} \\
&\leq \CC 2 \|\chi_m\|_{W^{1,\infty}(\Omega)} \sum_{x\in Z} \|R_m^{V,\pm} (I_X-\poop_{X_n}) P_V u\|_{H(\curl; B_{\CR  h_n }(x)\cap \Omega_\mp)} \\
&\leq \CC 2 \|\chi_m\|_{W^{1,\infty}(\Omega)} C_1 \|R_m^{V,\pm} (I_X-\poop_{X_n}) P_V u\|_{H(\curl; \tilde B)} \\
&\leq \CC 2 \|\chi_m\|_{W^{1,\infty}(\Omega)} C_1 \|R_m^{V,\pm}\|_{L(X_\pm,X_\mp)} \|(I_X-\poop_{X_n}) P_V \|_{L(X_n,X)} \|u\|_X.
\end{align*}
At last
\begin{align*}
\eqref{eq:termi} &=
\CC \sum_{x\in Z} \| (\chi_m-\chi_m(x)) (v-v_x)\|_{L^2(B_{\CR  h_n }(x)\cap \Omega_\mp)} \\
&+\CC \sum_{x\in Z} \| \Curl \chi_m \cdot (v-v_x)\|_{L^2(B_{\CR  h_n }(x)\cap \Omega_\mp)} \\
&+ \CC \sum_{x\in Z} \| (\chi_m-\chi_m(x)) \curl v\|_{L^2(B_{\CR  h_n }(x)\cap \Omega_\mp)} \\
&\leq \CC 2 \|\chi_m\|_{W^{1,\infty}(\Omega)} \sum_{x\in Z} \|(v-v_x)\|_{L^2(B_{\CR  h_n }(x)\cap \Omega_\mp)} \\
&+\CC \CR  h_n  \|\chi_m\|_{W^{1,\infty}(\Omega)} \sum_{x\in Z} \|\curl v\|_{L^2(B_{\CR  h_n }(x)\cap \Omega_\mp)} \\
&\leq \CC 2 \|\chi_m\|_{W^{1,\infty}(\Omega)} C_2 \hn^t \sum_{x\in Z} \|v\|_{H^t(B_{\CR  h_n }(x)\cap \Omega_\mp)} \\
&+\CC \CR  h_n  \|\chi_m\|_{W^{1,\infty}(\Omega)} \sum_{x\in Z} \|\curl v\|_{L^2(B_{\CR  h_n }(x)\cap\Omega_\mp)} \\
&\leq \CC 2 \|\chi_m\|_{W^{1,\infty}(\Omega)} C_2 \hn^t C_1 \|v\|_{H^t(\tilde B)} \\
&+\CC \CR  h_n  \|\chi_m\|_{W^{1,\infty}(\Omega)} C_1 \|\curl v\|_{L^2(\tilde B)} \\
&\lesssim h_n ^t \|u\|_X.
\end{align*}
Altogether we obtain
\begin{align*}
\| (T^{V,\pm}-T_n^{V,\pm})  &\poop_{X_n} P_V \|_{L(X_n,X)}
\lesssim h_n ^t + \|(I_X-\poop_{X_n}) P_V \|_{L(X_n,X)}
\end{align*}
which proves the claim.
\end{proof}

\begin{lemma}\label{lem:Wipop}
Let Assumptions~\ref{ass:abstractCutoff}, \ref{ass:abstractReflection} and \ref{ass:proj}, \ref{ass:localpol}, \ref{ass:localconst}, \ref{ass:Rconform} be satisfied.
Then for each $\signw\in\{0,-,+\}$ it holds
$\lim_{n\to\infty} \| (I_{W_0}-\ipop_{W_{0,n}}) T^{W_0,\signw} \poop_{W_{0,n}} P_{W_0} \|_n=0$.
\end{lemma}
\begin{proof}
Proceed as in the proof of Lemma~\ref{lem:Vipop}.
See also \cite[Theorem 5.5]{Halla:20PML}.
\end{proof}

\begin{theorem}\label{thm:Tcomp}
Let Assumptions~\ref{ass:abstractCutoff}, \ref{ass:abstractReflection} and \ref{ass:proj}, \ref{ass:localpol}, \ref{ass:localconst}, \ref{ass:Rconform} be satisfied.
Let $T(\cdot)$ and $T_n(\cdot)$ be as defined in \eqref{eq:DefT-identity}, \eqref{eq:DefT-classical}, \eqref{eq:DefT-crit} and \eqref{eq:DefTn} respectively.
Then for each $\lambda\in\Lambda_{R_V,R_{W_0}}$ it holds
$\lim_{n\to\infty} \|T(\lambda)-T_n(\lambda)\|_{L(X_n,X)}=0$.
\end{theorem}
\begin{proof}
Follows from $\lim_{n\to\infty} \|(I_X-P_{X_n}) K\|_{L(X)}=0$ for each compact operator $K\in L(X)$, Lemmata~\ref{lem:piPV}, \ref{lem:Vpoop}, \ref{lem:Wpoop}, \ref{lem:Vipop}, \ref{lem:Wipop} and the triangle inequality.
\end{proof}

\begin{remark}
Note that \cite{BonnetBDCarvalhoCiarlet:18} considers only lowest order finite elements and estimates the term similar to $(1-\ipop_{X_{\mp,n}})\chi_m R_m^{V,\pm}\poop_{X_n}P_Vu$ in a direct fashion.
Instead, one can apply the discrete commutator technique \cite{Bertoluzza:99} as in Lemma \ref{lem:Vipop} to generalize the results of \cite{BonnetBDCarvalhoCiarlet:18} to high order methods.
\end{remark}

\subsection{Explicit construction}\label{subsec:Tcomp-explicit}

\begin{definition}[locally $R$-conform mesh (Definition 3 of \cite{BonnetBDCarvalhoCiarlet:18})]
\label{def:locallyRconformmesh}
For a mesh $\calT$ of $\Omega$ let $\calT^\pm_n:=\{\Delta\in\calT\colon \Delta\cap \Omega_n\cap\Omega_\pm\neq\emptyset\}$, $n=1,\dots,2N$.
A mesh $\calT$ is called locally $R$-conform, if for all $\Delta\in\calT$ either $\Delta\subset\Omega_+$ or $\Delta\subset\Omega_-$, and if the image of each $\Delta\in\calT^\pm_n$ by the underlying the geometrical transformations of $R^{W_0,\pm}_n$
($\phi_{n,\pm,k,m}$, $k=1,\dots,p_{n,\mp}$, $m=1,\dots,p_{n,\pm}$, $n=1,\dots,N$ and $\phi_{n,\pm}$, $n=N+1,\dots,2N$)
belongs to $\calT^\mp_n$, $n=1,\dots,2N$.
\end{definition}
For the construction of such meshes we refer to \cite{BonnetBDCarvalhoCiarlet:18} and Figures~\ref{fig:geometry} and \ref{fig:coarsemesh}.
We note that for a technical reason our Definition~\ref{def:locallyRconformmesh} is slightly different to \cite[Definition 3]{BonnetBDCarvalhoCiarlet:18}.
We formulate our definition in terms of the slightly lager patches $\Omega_n$ instead of $S_n$.

\begin{lemma}\label{lem:discretespaces}
Assume that $\Omega$ is a Lipschitz polyhedron and that all angles $(\alpha_n)_{n=1,\dots,N}$ of the interface $\Sigma$ are in $2\pi\setQ$.
Let $\chi_n$ and $R_n^{V,\pm}, R_n^{W_0,\pm}$ be as defined in Subsection~\ref{subsec:Tcomp-explicit}.
Let $(\calT_n)_{n\in\setN}$ be a sequence of shape-regular and locally $R$-conform triangular meshes of $\Omega$ with maximal element diameter $h_n \to0$ as $n\to\infty$.
Let $X_n$ be a N\'ed\'elec finite element space on the mesh $\calT_n$ such that the polynomial degree $k\in\setN$ is uniform on all triangles $\Delta\in\bigcup_{n=1}^{2N}\Omega_n$.
Then Assumptions~\ref{ass:proj}, \ref{ass:localpol}, \ref{ass:localconst} and \ref{ass:Rconform} are satisfied.
\end{lemma}
\begin{proof}
Let $(W_{0,n}, X_n, J_n)$ be the subcomplex of $(W_0,X,J)$ associated to $X_n$, see e.g.\ \cite{Zaglmayr:06}.
Then the existence of $L^2$-uniformly bounded cochain projections $\poop_{W_{0,n}}$, $\poop_{X_n}$, $\poop_{J_n}$ follows as in \cite[Theorem 5.9]{ArnoldFalkWinther:10}.
Hence Assumption~\ref{ass:proj} is satisfied.
We can use the Scott-Zhang-interpolant as local projection $\ipop_{W_{0,n}}$, see e.g.\ \cite[Lemma 1.130]{ErnGuermond:04}.
We can use the standard $H(\curl)$-interpolant as local projection $\ipop_{X_n}$, which is defined in terms of the degrees of freedom and basis functions \cite[p.\ 62, 80]{Zaglmayr:06}.
Here we exploit that in the special $2D$ setting the $H(\curl)$-interpolant is bounded on $H(\curl)$.
Indeed we recall that the degrees of freedom are
\begin{align*}
\spl q_l, \tv\cdot u \spr_{H^{1/2}(e)\times H^{-1/2}(e)}, \quad
(q_l)_{l=0,\dots,k}\text{ basis of } P^k(e),
\end{align*}
for each edge $e$ in the skeleton of $\calT_n$, and
\begin{align*}
\spl s_l, \curl u \spr_{L^2(\Delta)}, \quad
s_l &\text{ basis for }P^{k-1}(\Delta)/\setR,\\
\spl \nabla p_l, u \spr_{L^2(\Delta)}, \quad
p_l &\text{ basis for }\lambda_1\lambda_2\lambda_3 P^{k-2}(\Delta),
\end{align*}
for each triangle $\Delta\in\calT_n$.
Hereby $P^k$ denotes the space of polynomials of degree lower equal than $k$ and $\lambda_j$, $j=1,2,3$ denote the barycentric coordinates.
Thus Assumption~\ref{ass:localpol} is satisfied.
It can easily be seen that Assumption~\ref{ass:localconst} is satisfied by finite element spaces.
Since the meshes $\calT_n$ are locally $R$-conform and the polynomial degree $k$ is locally uniform too, Assumption~\ref{ass:Rconform} is satisfied as well.
\end{proof}

\subsection{Summary of results}\label{subsec:Tcomp-results}

\begin{theorem}
Assume that $\Omega$ is a Lipschitz polyhedron and that all angles $(\alpha_n)_{n=1,\dots,N}$ of the interface $\Sigma$ are in $2\pi\setQ$.
Let $\lambda\in\Lambda_\alpha$ be in the resolvent of $A(\cdot)$, $f\in X$ and $u\in X$ be the solution to $A(\lambda)u=f$.

Let $(\calT_n)_{n\in\setN}$ be a sequence of shape-regular and locally $R$-conform triangular meshes of $\Omega$ with maximal element diameter $h_n \to0$ as $n\to\infty$.
Further let $k_0\in\setN$ and $X_n$ be a N\'ed\'elec finite element space on the mesh $\calT_n$ with maximal local polynomial degrees $k\leq k_0$ and such that the local polynomial degree $k\in\setN$ is uniform on all triangles $\Delta\cap\bigcup_{n=1}^{2N}\Omega_n\neq\emptyset$.

Then there exist constants $C,n_0>0$ such that for all $n>n_0$ there exists a unique solution $u_n\in X_n$ to the approximated equation $P_{X_n}A(\lambda)u_n=P_{X_n}f$ which satisfies
$\|u-u_n\|_{H(\curl;\Omega)} \leq C \inf_{u_n'\in X_n} \|u-u_n'\|_{H(\curl;\Omega)}$.
\end{theorem}
\begin{proof}
Follows from Theorems \ref{thm:wTc-explicit}, \ref{thm:Tcomp}, Lemma \ref{lem:discretespaces}, \cite[Corollary 2.8, iii)]{Halla:19Tcomp} and basic approximation theory.
\end{proof}

\begin{theorem}\label{thm:MainResult}
Assume that $\Omega$ is a Lipschitz polyhedron and that all angles $(\alpha_n)_{n=1,\dots,N}$ of the interface $\Sigma$ are in $2\pi\setQ$.
Consider the eigenvalue problem to find $(\lambda,u)\in\Lambda_\alpha\times X\setminus\{0\}$ such that $A(\lambda)u=0$.

Let $(\calT_n)_{n\in\setN}$ be a sequence of shape-regular and locally $R$-conform triangular meshes of $\Omega$ with maximal element diameter $h_n \to0$ as $n\to\infty$.
Further let $k_0\in\setN$ and $X_n$ be a N\'ed\'elec finite element space on the mesh $\calT_n$ with maximal local polynomial degrees $k\leq k_0$ and such that the local polynomial degree $k\in\setN$ is uniform on all triangles $\Delta\cap\bigcup_{n=1}^{2N}\Omega_n\neq\emptyset$.
Consider the discrete eigenvalue problem to find $(\lambda_n,u_n)\in\Lambda_\alpha\times X_n\setminus\{0\}$ such that $P_{X_n}A(\lambda_n)u_n=0$.

Then the eigenvalues $\lambda_n$ and eigenfunctions $u_n$ converge to $\lambda$ and $u$ respective in the sense of \cite[i)-vii)]{Halla:19Tcomp}.
For example, if $\tilde\Lambda\subset\Lambda_\alpha$ is a compact set with rectifiable boundary
$\partial\tilde\Lambda\subset\rho\big(A(\cdot)\big)$ and a simple eigenvalue $\{\lambda\}=\tilde\Lambda\cap\sigma\big(A(\cdot)\big)$ with normalized eigenfunction $u\in X$, then there exist $n_0\in\setN$ and $c>0$ such that for all $n>n_0$ the set $\tilde\Lambda\cap\sigma\big(P_{X_n}A(\cdot)|_{X_n}\big)$ consists of a simple eigenvalue $\lambda_n$ with normalized eigenfunction $u_n\in X_n$ and
\begin{align*}
|\lambda-\lambda_n|\leq c \inf_{u'_n\in X_n} \|u-u'_n\|_X^2,
\quad
\inf_{\mu\in\setC,|\mu|=1}\|\mu u-u_n\|_X\leq c \inf_{u'_n\in X_n} \|u-u'_n\|_X.
\end{align*}
\end{theorem}
\begin{proof}
Follows from Theorems \ref{thm:wTc-explicit}, \ref{thm:Tcomp}, Lemma \ref{lem:discretespaces} and \cite[Corollary 2.8]{Halla:19Tcomp}.
\end{proof}

To compute the eigenvalues of the rational matrix function $A_n(\cdot)$ it is most convenient to linearize it and subsequently to apply a linear eigenvalue solver.
To this end let
\begin{align*}
Y_n:=\{(\curl u)|_{\Omega_-}\colon u\in X_n\} \subset L^2(\Omega_-)
\end{align*}
with associated orthogonal projection $p_{Y_n}\in L(L^2(\Omega_-),Y_n)$ and
\begin{align*}
v:=\frac{\omega_\mu}{\sqrt{\mu_-}}\frac{1}{\lambda-\omega_\mu^2} \curl u|_{\Omega_-}.
\end{align*}
Then for $\lambda\neq\omega_\mu^2$ the function $u\in X_n$, solves $A_n(\lambda)u=0$ if and only if $(u,v)\in X_n\times Y_n$ solves
\begin{align}\label{eq:LEVP}
\begin{aligned}
0&=\spl \mu_+^{-1} \curl u,\curl u' \spr_{L^2(\Omega_+)}
+\spl \mu_-^{-1} \curl u,\curl u' \spr_{L^2(\Omega_-)}\\
&-\lambda \spl \epsilon_+ \curl u,\curl u' \spr_{L^2(\Omega_+)}
+\omega_\epsilon^2 \spl \epsilon_- \curl u,\curl u' \spr_{L^2(\Omega_-)}\\
&+\omega_\mu \spl \mu_-^{-1/2} v,\curl u' \spr_{L^2(\Omega_-)}
+\omega_\mu \spl \mu_-^{-1/2} \curl u,v' \spr_{L^2(\Omega_-)}
-(\lambda-\omega_\mu^2)\spl v,v' \spr_{L^2(\Omega_-)}
\end{aligned}
\end{align}
for all $(u',v')\in X_n\times Y_n$.
However, if $\mu_-$ is not constant then the above correspondence is not true.
In this case one can choose e.g.\ $v:=\frac{1}{\lambda-\omega_\mu^2} \curl u|_{\Omega_-}$ as auxiliary variable.
Then the correspondence still holds, but the obtained matrix stencil is not selfadjoint.
Alternatively, we can analyze the selfadjoint matrix stencil in \eqref{eq:LEVP} as follows.
If $(u,v)\in X_n\times Y_n$ solves \eqref{eq:LEVP}, then $u\in X_n$ solves $(A_n(\lambda)+K_n(\lambda))u=0$ with $K_n(\lambda)\in L(X_n)$ defined by
\begin{align*}
\spl K_n(\lambda)u,u'\spr_X:=
\frac{\omega_\mu^2}{\lambda-\omega_\mu^2}
\spl \mu_-^{-1/2} (1-p_{Y_n}) (\mu_-^{-1/2} \curl u), \curl u' \spr_{L^2(\Omega_-)}
\end{align*}
for all $u,u'\in X_n$.
For simplicity let $X_n$ have a uniform polynomial degree $k\geq1$.
We compute
\begin{align*}
\spl \mu_-^{-1/2} (1-p_{Y_n}) &(\mu_-^{-1/2} \curl u), \curl u' \spr_{L^2(\Omega_-)}\\
&=\sum_{\Delta\in\calT_n, \Delta\subset\Omega_-}
\spl \mu_-^{-1/2} (1-p^{k-1}_\Delta) (\mu_-^{-1/2} \curl u)|_\Delta, \curl u' \spr_{L^2(\Delta)}
\end{align*}
with orthogonal projections $p^{k-1}_\Delta\in L(L^2(\Delta),P^{k-1}(\Delta))$.
Since for any constant $c_\Delta\in\setC$
\begin{align*}
(1-p^{k-1}_\Delta) (\mu_-^{-1/2} \curl u)|_\Delta=(1-p^{k-1}_\Delta) (\mu_-^{-1/2}-c_\Delta) \curl u)|_\Delta \quad\text{for }u\in X_n
\end{align*}
we can estimate
\begin{align*}
|\spl \mu_-^{-1/2} (1-p^{k-1}_\Delta) (\mu_-^{-1/2} &\curl u)|_\Delta, \curl u' \spr_{L^2(\Delta)}|\\
&\leq \|\mu^{-1/2}\|_{L^\infty(\Omega_-)} \|\mu_-^{-1/2}-c_\Delta\|_{L^\infty(\Delta)} \|\curl u\|^2_{L^2(\Delta)}
\end{align*}
\sloppy
If there exists $s\in(0,1]$ such that $\mu_-\in W^{s,\infty}(\Omega_-)$ or equivalently $\mu^{-1/2}_-\in W^{s,\infty}(\Omega_-)$, then we can choose $c_\Delta$ such that
$\|\mu^{-1/2}-c_\Delta\|_{L^\infty(\Delta)}$ $\leq$ $h_n^s \|\mu^{-1/2}\|_{W^{s,\infty}(\Omega_-)}$.
Thus
\begin{align*}
\|K_n(\lambda)\|_{L(X_n)} \leq h_n^s \|\mu^{-1/2}\|_{W^{s,\infty}(\Omega_-)} \|\mu^{-1/2}\|_{L^\infty(\Omega_-)} \frac{\omega_\mu^2}{|\lambda-\omega_\mu^2|}
\end{align*}
and hence $\lim_{n\to\infty} \|K_n(\lambda)\|_{L(X_n)}=0$.

\section{Computational experiments}\label{sec:compexamp}

\begin{figure}
\begin{subfigure}{.47\textwidth}
\centering
\includegraphics[clip, trim={15cm 2cm 15cm 3cm}, scale=0.15]{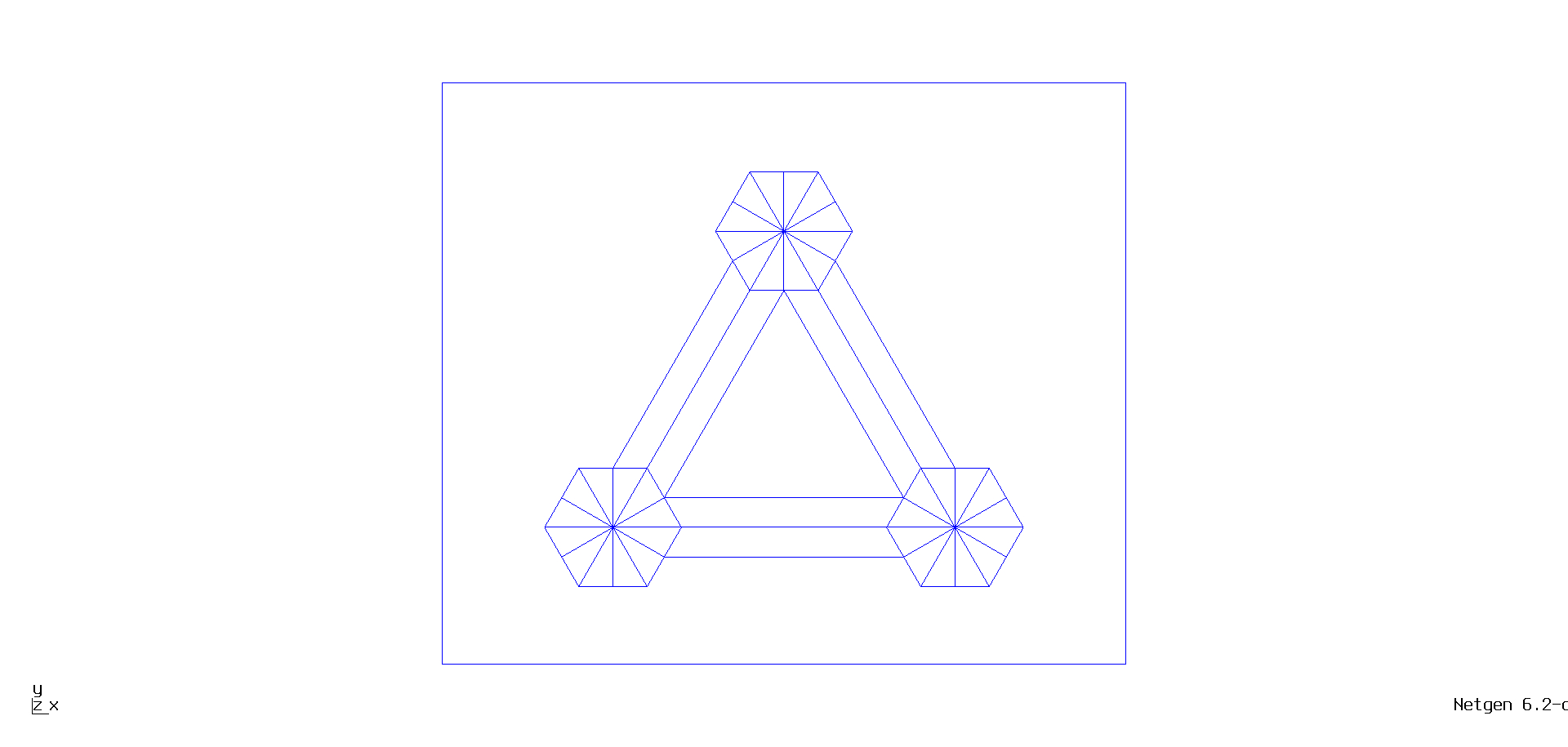}
\caption{Geometrical setup to build locally $R$-conform meshes.}
\label{fig:geometry}
\end{subfigure}
\hfill
\begin{subfigure}{.47\textwidth}
\centering
\includegraphics[clip, trim={15cm 1.5cm 16cm 4cm}, scale=0.15]{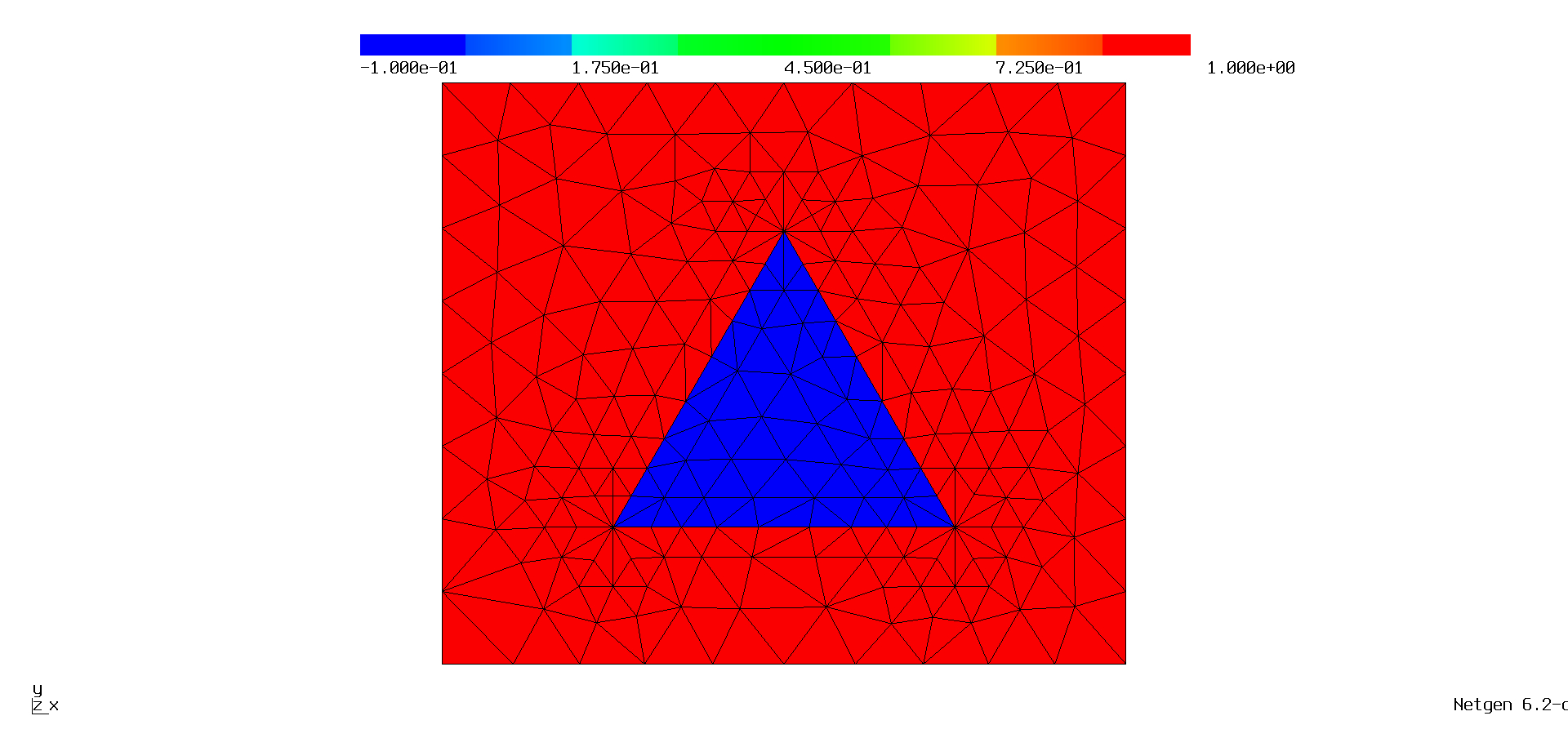}
\caption{Coarse locally $R$-conform mesh. The color coding highlights $\Omega_\pm$.}
\label{fig:coarsemesh}
\end{subfigure}
\caption{Construction of locally $R$-conform meshes.}
\end{figure}

\sloppy
For our forthcoming computational experiments we consider a rectangular domain $\Omega=(-0.5,1.5)$ $\times$ $(-0.5,1.3)$ with $\Omega_-$ being the equilateral triangle with the corners $(0,0)$, $(1,0)$ and $(\cos(\pi/3), \sin(\pi/3))$, see Figure~\ref{fig:coarsemesh}.
Thus we obtain the critical interval $I_\alpha=[-5,-1/5]$.
In order to build locally $R$-conform meshes we proceed as in \cite{BonnetBDCarvalhoCiarlet:18}.
For each corner $c_n$, $n=1,2,3$ we construct a neighborhood $\Omega_n$ of $c_n$ as a regular convex six-sided polygon with center $c_n$ and edges aligned to the edges of $\Omega_-$.
In addition we symmetrically bisect each of the six slices of $\Omega_n$.
Then for each edge $e_n$, $n=1,2,3$ we construct two symmetric trapezoids which have corner points of the neighboring $\partial\Omega_m$ as in Figure~\ref{fig:geometry}.
For our computational experiments we use the mesh generator Netgen and the finite element software NG-Solve.

\subsection{Source problem}

\begin{figure}
\begin{subfigure}{.48\textwidth}
\centering
\includegraphics[clip, trim={15cm 0cm 10cm 0cm}, scale=0.14]{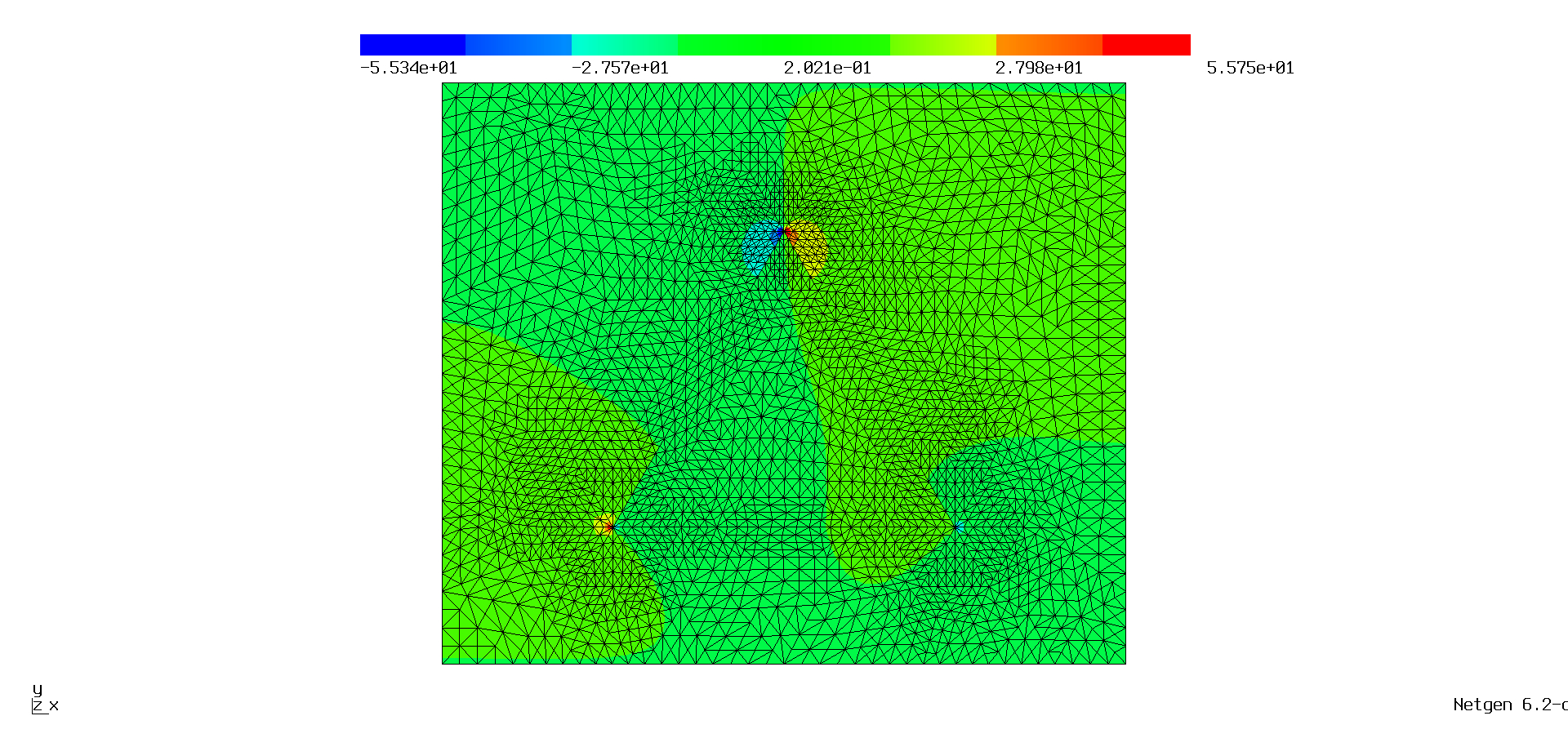}
\caption{First component of $u_n$.}
\label{fig:u1}
\end{subfigure}
\hfill
\begin{subfigure}{.48\textwidth}
\centering
\includegraphics[clip, trim={15cm 0cm 10cm 0cm}, scale=0.14]{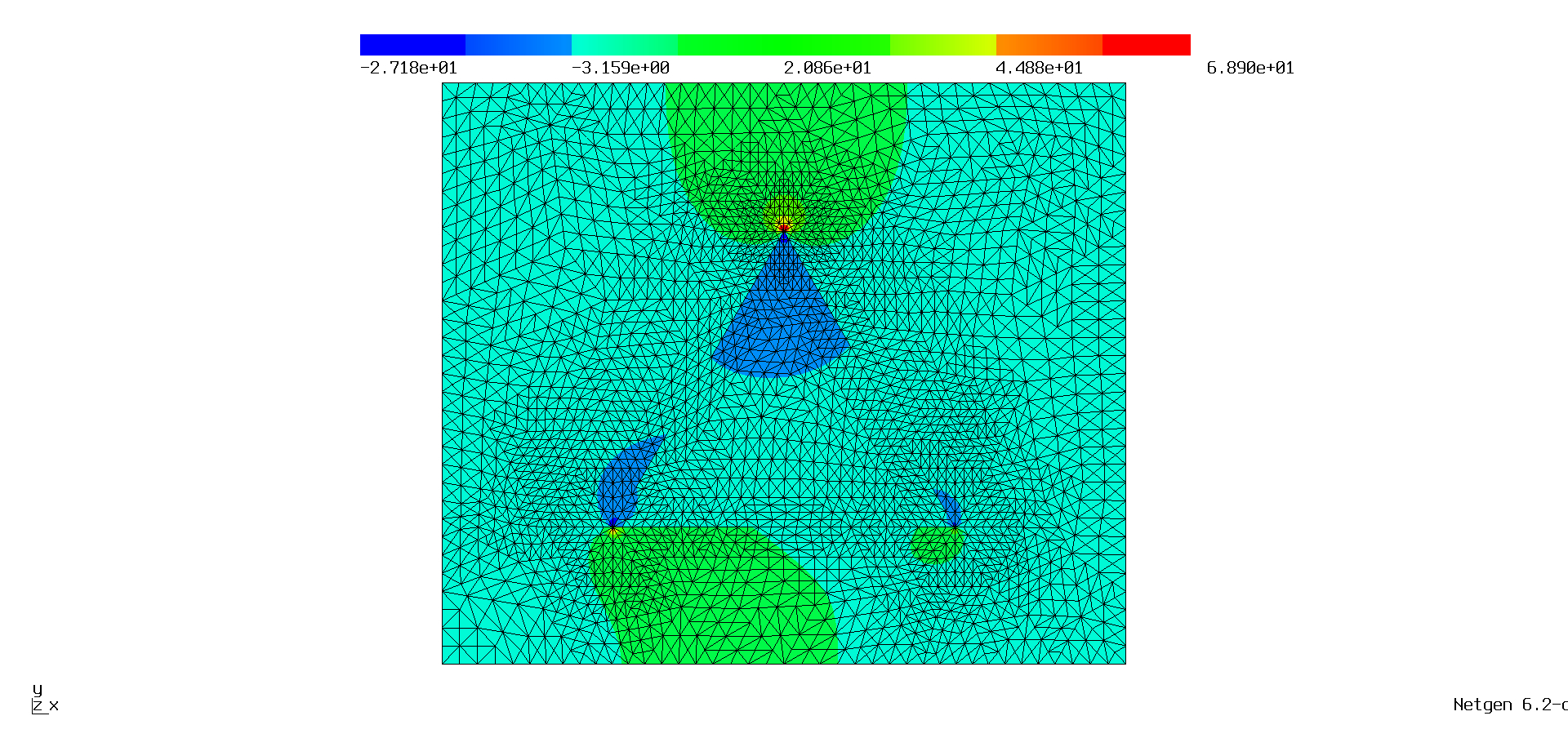}
\caption{Second component of $u_n$.}
\label{fig:u2}
\end{subfigure}
\caption{Computed solution $u_n$ for $f=(1,1)^\top$.}
\label{fig:u}
\end{figure}

We consider the source problem \eqref{eq:PDE} with right hand side $f=(1,1)^\top$, $\lambda=1$ and $\omega_\mu=\omega_\epsilon=\sqrt{2}$, $\mu_+=1$, $\mu_-=1/10$, $\epsilon_+=1$, $\epsilon_-=10$.
Thus we obtain the contrasts $\kappa_{\mu^{-1}}(\lambda)=\kappa_\epsilon(\lambda)=-10$.
We consider a sequence of locally $R$-conform meshes $(\calT_n)_{n\setN}$ obtained from uniform refinements of an initial mesh $\calT_0$ with maximal element diameter $h_{\max}=0.2$.
Subsequently we build $X_n$ as lowest order N\'ed\'elec elements on $\calT_n$.
As we do not have an analytical solution at our disposal, we take the finite element solution on the finest mesh as reference.
In Figure~\ref{fig:errors} we present the relative errors in the $H(\curl)$ and $L^2$ norms.
To gain more confidence in our approximations we compute a second reference solution in an alternative way.
By means of a Helmholtz decomposition we can represent $u=\epsilon(\lambda)^{-1}\Curl v$ with $v\in H^1(\Omega)$ being the solution to
\begin{subequations}\label{eq:PDE-pot}
\begin{align}
-\div(\epsilon(\lambda)^{-1}\nabla v)-\lambda \mu(\lambda)v &= \mu(\lambda) f_0 \quad\text{in }\Omega,\\
\nv\cdot\epsilon(\lambda)^{-1}\nabla v&=0 \quad\text{on }\partial\Omega,
\end{align}
\end{subequations}
with potential $f_0(x)=x_1-x_2$, $\Curl f_0=f$.
Recall that the approximation of \eqref{eq:PDE-pot} was already analyzed in \cite{BonnetBDCarvalhoCiarlet:18}.
We compute an approximate solution $v_n\in H^1(\Omega)$ to \eqref{eq:PDE-pot} with scalar finite elements on the finest mesh $\calT_n$ and set $u_n':=\epsilon(\lambda)^{-1}\Curl v_n \in (L^2(\Omega))^2$.
Consequently we can compare $u_n$ to $u_n'$ in the $L^2$ norm, see the last column in Figure~\ref{fig:errors}.
In accordance with our theory we observe decreasing errors in Figure~\ref{fig:errors}.
In Figure~\ref{fig:u} we plot the computed solution.

\begin{figure}
\centering
\begin{tabular}{c|c|c|c|c}
n & \# dof & $\frac{\|u_6-u_n\|_{H(\curl)}}{\|u_6\|_{H(\curl)}}$ & $\frac{\|u_6-u_n\|_{L^2}}{\|u_6\|_{L^2}}$ & $\frac{\|u'_6-u_n\|_{L^2}}{\|u'_6\|_{L^2}}$ \\
\hline
0 & 1,278 & 4.793e-01 & 5.783e-01 & 5.908e-01 \\
1 & 5,688 & 2.683e-01 & 3.346e-01 & 3.433e-01 \\
2 & 23,202 & 1.655e-01 & 2.120e-01 & 2.185e-01 \\
3 & 93,006 & 1.092e-01 & 1.427e-01 & 1.481e-01 \\
4 & 371,718 & 7.485e-02 & 9.900e-02 & 1.042e-01 \\
5 & 1,485,558 & 5.184e-02 & 6.901e-02 & 7.532e-02 \\
6 & 5,938,902 & 0.0 & 0.0 & 3.305e-02
\end{tabular}
\caption{Computed errors.}
\label{fig:errors}
\end{figure}

\subsection{Eigenvalue problem}

\begin{figure}
\centering
\includegraphics[clip, trim={5cm 0cm 0cm 0cm}, scale=0.3]{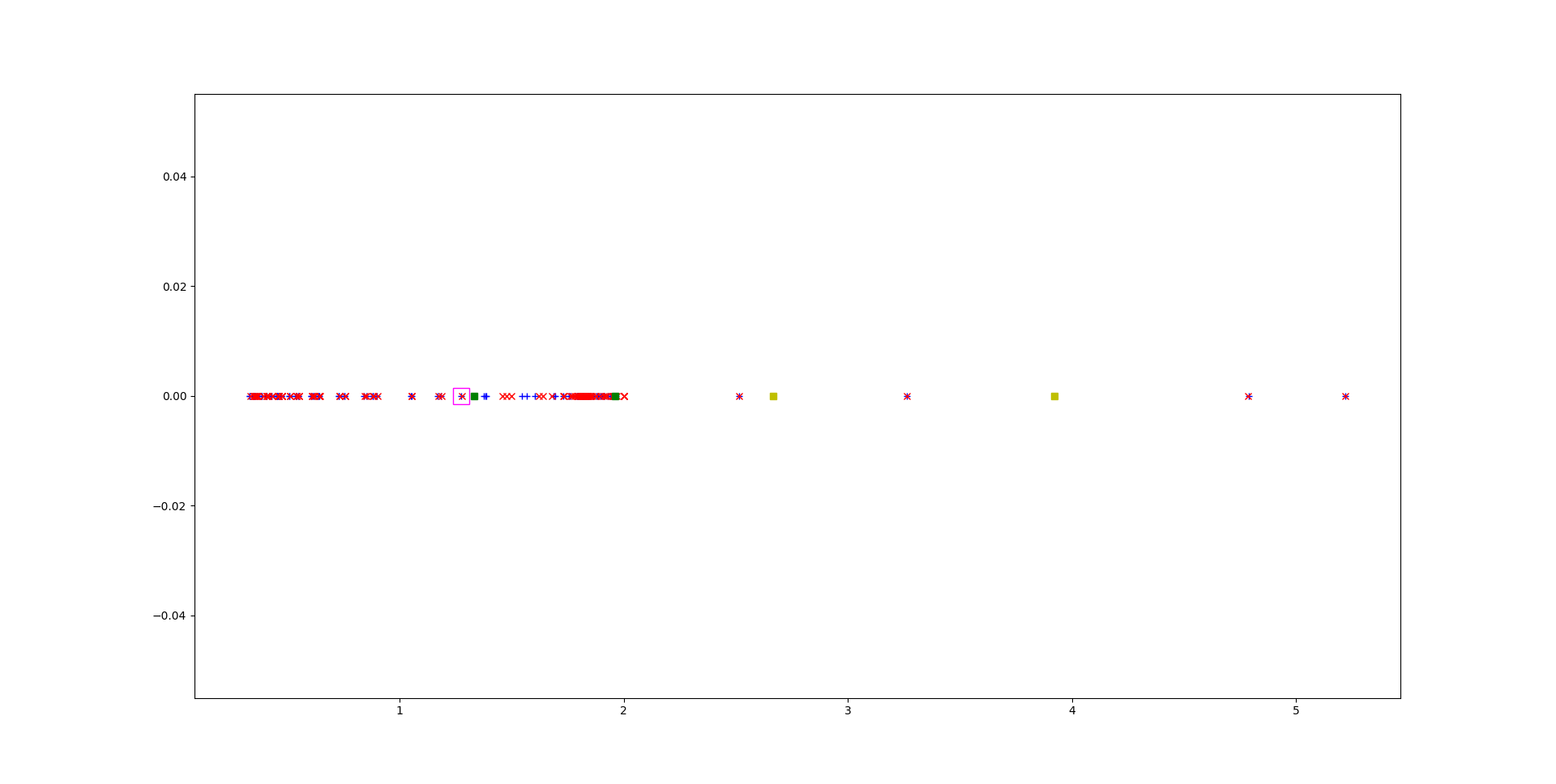}
\caption{Computed spectrum. Blue + mark the $H(\curl)$ computations. Red x mark the $H^1$ computations. Colored squares mark $4/3$, $100/51$ and $8/3, 200/51$. Magenta box marks the biggest eigenvalue below $4/3$.}
\label{fig:spec}
\end{figure}

\begin{figure}
\begin{subfigure}{.48\textwidth}
\centering
\includegraphics[clip, trim={15cm 0cm 10cm 0cm}, scale=0.14]{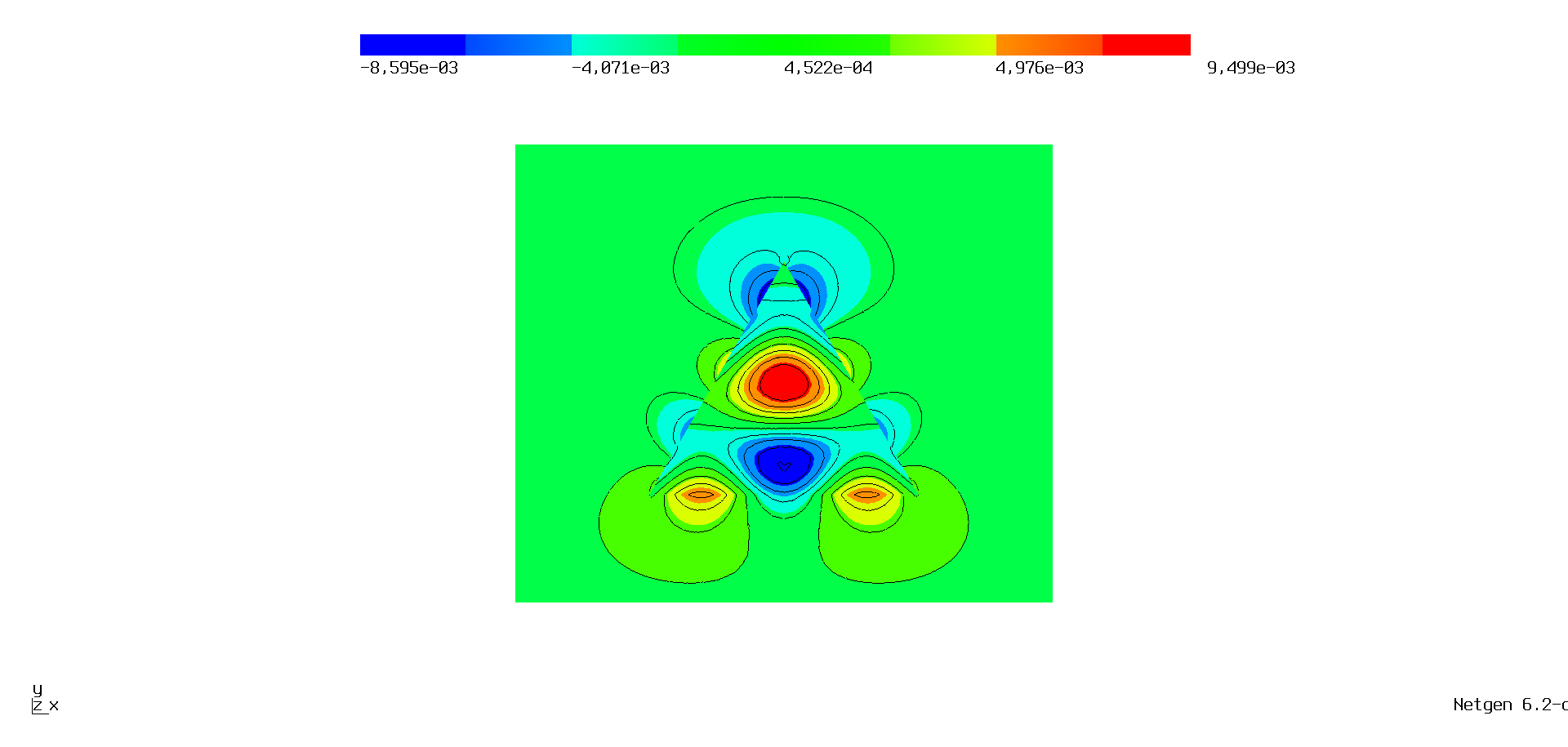}
\caption{First component of $u_n$.}
\label{fig:ef1}
\end{subfigure}
\hfill
\begin{subfigure}{.48\textwidth}
\centering
\includegraphics[clip, trim={15cm 0cm 10cm 0cm}, scale=0.14]{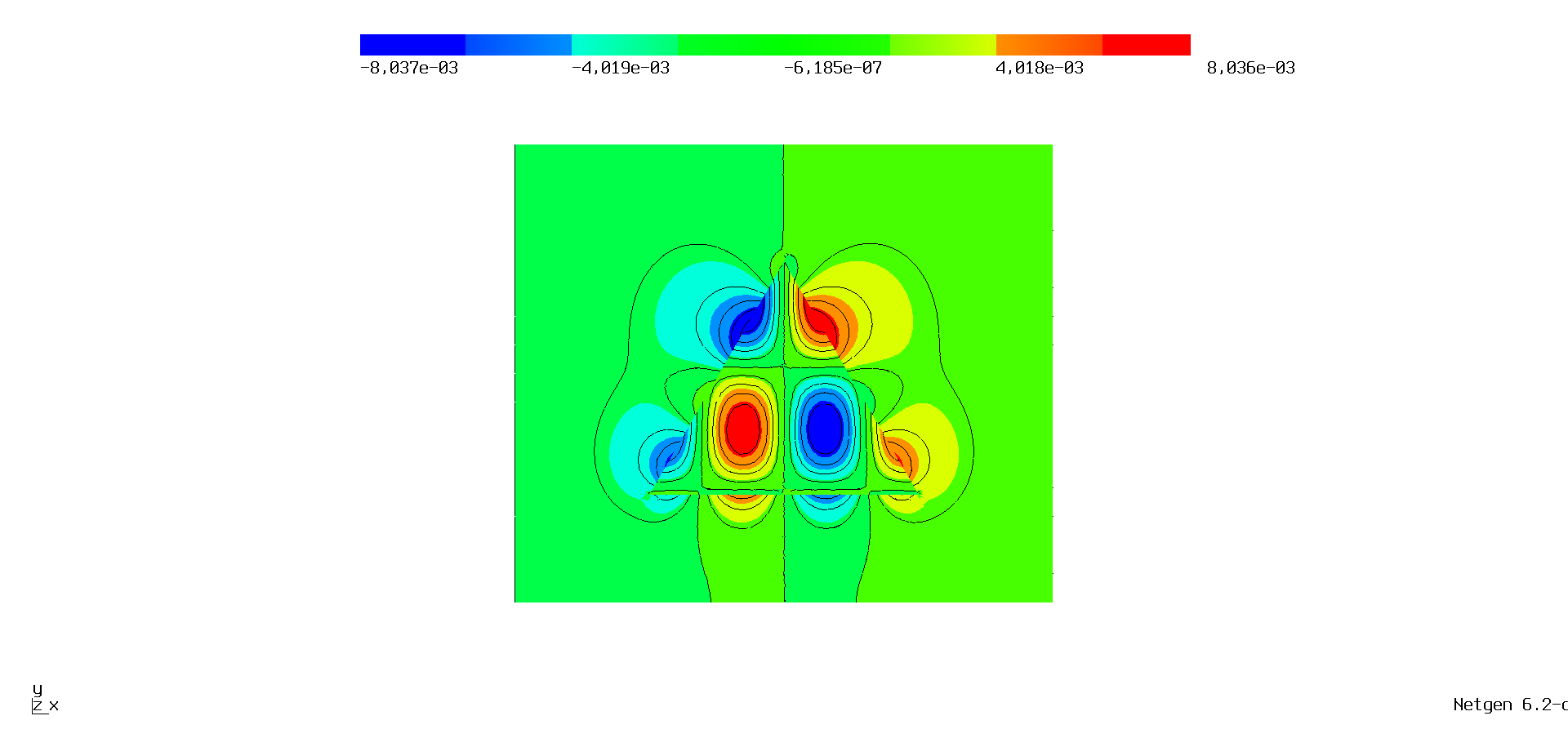}
\caption{Second component of $u_n$.}
\label{fig:ef2}
\end{subfigure}
\caption{Computed eigenfunction $u_n$ (n=2) to the eigenvalue $\lambda\approx 1.276$.}
\label{fig:ef}
\end{figure}

Let $\omega_\mu=2$, $\omega_\epsilon=\sqrt{2}$, $\mu_+=1$, $\mu_-=10$, $\epsilon_+=1$, $\epsilon_-=10$.
Then $\kappa_{\mu^{-1}}(\lambda)\in I_\alpha$ if and only if $\lambda\in [8/3,200/51]$ and
$\kappa_{\epsilon}(\lambda) \in I_\alpha$ if and only if $\lambda\in [4/3,100/51]$.
We consider the eigenvalue problem \eqref{eq:PDE}, i.e. we search for $(\lambda,u)\in\setC\times X\setminus\{0\}$ which solve \eqref{eq:PDE} with $f=0$.
Since this is a rational eigenvalue problem, we linearize it with the auxiliary variable $v:=\frac{\omega_\mu}{\sqrt{\mu_-}}\frac{1}{\lambda-\omega_\mu^2} \curl u|_{\Omega_-}$ as in Section~\ref{subsec:Tcomp-results}.
Subsequently we apply the Arnoldi algorithm of NG-Solve to solve the linear matrix eigenvalue problem.
As for the source problem we use the associated equation for the $\Curl$-potential to compute a second reference solution.
In Figure~\ref{fig:spec} we present the computed spectrum.
Outside the interval $[4/3,100/51]$ we note a good accordance of the eigenvalues between the two computations.
Inside the interval $[4/3,100/51]$ (between the green boxes) we note an accumulation of eigenvalues, which is to be expected since the operator function is not Fredholm in this interval.
In addition we also observe a good accordance of the eigenvalues inside $[8/3,200/51]$, which is not covered by our theory.
To compute errors we choose the biggest eigenvalue below $4/3$ (marked with a magenta box) and list in Figure~\ref{fig:errors-ev} the associated errors.
In accordance with our theory we observe decreasing errors in Figure~\ref{fig:errors-ev}.
In Figure~\ref{fig:ef} we plot the corresponding eigenfunction.

\begin{figure}
\centering
\begin{tabular}{c|c|c|c}
n & \# dof & $\frac{|\lambda_6-\lambda_n|}{|\lambda_6|}$ & $\frac{|\lambda'_6-\lambda_n|}{|\lambda'_6|}$ \\
\hline
0 & 1,692 & 3.728e-02 & 3.729e-02 \\
1 & 7,344 & 1.185e-02 & 1.187e-02 \\
2 & 29,826 & 3.066e-03 & 3.078e-03 \\
3 & 119,502 & 7.633e-04 & 7.755e-04 \\
4 & 477,702 & 1.820e-04 & 1.942e-04 \\
5 & 1,909,494 & 3.643e-05 & 4.858e-05 \\
6 & 7,634,646 & 0.0 & 1.215e-05 
\end{tabular}
\caption{Computed errors in $\lambda\approx 1.276$.} 
\label{fig:errors-ev}
\end{figure}

\section{Conclusion}\label{sec:conclusion}

In this article we extended the analysis of \cite{BonnetBDCarvalhoCiarlet:18} to dispersive time-harmonic two-dimensional vectorial electromagnetic wave equations.
For the definition of our $T$-operator we introduced reflection operators for each Helmholtz component, whereby both individual reflection operators are based on the geometric construction in \cite{BonnetBDCarvalhoCiarlet:18}.
We showed that h-finite element methods with N\'ed\'elec elements and locally $R$-conform meshes are $T$-compatible.
Consequently the $T$-compatibility ensures \cite{Halla:19Tcomp} that the holomorphic eigenvalue problems are approximated reliably.
We underlayed our results with computational experiments.
However, our analysis requires an additional assumption on the contrast $\kappa_{\mu^{-1}}$, which is not necessary for the results on the continuous level in 2D \cite{BonnetBDChesnelCiarlet:14b}.
Indeed, our computational experiments suggest that the finite element method yields reliable results even if the condition on $\kappa_{\mu^{-1}}$ is violated.
%

\bibliographystyle{amsplain}
\bibliography{short_biblio}

\end{document}